\PassOptionsToPackage{pdfusetitle, psdextra}{hyperref}
\PassOptionsToPackage{capitalize}{cleveref}
\documentclass{article}
\pdfoutput=1

\usepackage{subfiles}
\usepackage{hyperref}
\usepackage{biblatex}
\newcommand{\myCite}{\parencite}
\usepackage{amsthm}
\usepackage{amsmath}
\usepackage{cleveref}
\usepackage{my-preamble}

\newtheorem{myDummyTheorem}{Dummy}[section]
\theoremstyle{definition}
\newtheorem{myDefinition}[myDummyTheorem]{Definition}
\newtheorem{myConstruction}[myDummyTheorem]{Construction}
\newtheorem{myNotation}[myDummyTheorem]{Notation}
\newtheorem{myConvention}[myDummyTheorem]{Convention}
\theoremstyle{remark}

\newtheorem{myRemark}[myDummyTheorem]{Remark}
\theoremstyle{plain}
\newtheorem{myProposition}[myDummyTheorem]{Proposition}
\newtheorem{myLemma}[myDummyTheorem]{Lemma}
\newtheorem{myTheorem}[myDummyTheorem]{Theorem}
\newtheorem*{myTheorem*}{Theorem}
\newtheorem{myCorollary}[myDummyTheorem]{Corollary}

\newenvironment{myProof}{\begin{proof}}{\end{proof}}
\newcommand{\myQED}{\qed}

\newcommand{\myRef}{\cref}
\newcommand{\myLabelRef}{\labelcref}

\title{\myTitle}
\author{\myAuthor}
\date{October 23, 2024}

\begin{document}

\maketitle

\begin{abstract}

We propose a definition of higher inductive types in
$(\infty,1)$-categories with finite limits. We show that the
$(\infty,1)$-category of $(\infty,1)$-categories with higher inductive types is
finitarily presentable. In particular, the initial
$(\infty,1)$-category with higher inductive types exists. We prove a form
of canonicity: the global section functor for the initial
$(\infty,1)$-category with higher inductive types preserves higher
inductive types.

\end{abstract}

\section{Introduction}

\emph{Higher inductive types} are one of features of homotopy type
theory \myCite{hottbook}. An ordinary inductive type is a type freely
generated by a specified set of constructors for elements of the
type. An example of inductive type is the type of natural numbers,
whose constructors are zero and successor. A higher inductive type
may, in addition, have constructors for identifications between
elements of the type. Examples of higher inductive type are the
circle, whose constructors are the base point and the loop on it,
propositional truncations, where any two elements are forced to be
identified by a constructor, and set truncations, where any two
parallel identifications are forced to be identified by a
constructor. Higher inductive types allow us to construct various
types including: colimits such as \(n\)-spheres; disjunctions and
existential quantifications of propositions, which are combinations of
propositional truncations and colimits; homotopy groups, which are the
set truncations of loop spaces.

Homotopy type theory is expected to be an internal language for
\((\infty,1)\)-categories with some structures. For example, theories with
identity types are equivalent to \((\infty,1)\)-categories with finite
limits \myCite{kapulkin2018homotopy,nguyen2022type-arxiv}. A natural
question is: what \((\infty,1)\)-categorical structure corresponds to
higher inductive types? There have been some studies on semantics of
higher inductive types in specific contexts:
\myCite{lumsdaine2020semantics} in the context of model categories;
\myCite{coquand2018higher} in the cubical set model. In the present
paper, we propose a definition of higher inductive types in
\((\infty,1)\)-categories with finite limits.

Our definition of higher inductive types has the following features.
\begin{enumerate}
\item \label[feature]{property-hit-universal-property} Higher inductive types
  are defined by universal properties.
\item \label[feature]{property-hit-presentable} The \((\infty,1)\)-category of
  \((\infty,1)\)-categories with higher inductive types is finitarily
  presentable.
\end{enumerate}
These are summarized in
\myRef{prop-ind-type-pres}. \myRef{property-hit-universal-property}
implies that, like limits and colimits, higher inductive types are
unique if they exist. In particular, the \((\infty,1)\)-category of
\((\infty,1)\)-categories with higher inductive types is a
sub-\((\infty,1)\)-category of the \((\infty,1)\)-category of
\((\infty,1)\)-categories with finite limits. By
\myRef{property-hit-presentable}, the \((\infty,1)\)-category of
\((\infty,1)\)-categories with higher inductive types has small limits and
colimits and admits various free constructions. In particular, we have
the initial \((\infty,1)\)-category with higher inductive types which is
expected to be the syntactic \((\infty,1)\)-category of the type theory
with higher inductive types.

We prove a form of \emph{canonicity} for higher inductive
types. Canonicity in type theory means that every closed term of an
inductive type is identical to a canonical form. For example, the
canonicity for the type of natural numbers asserts that every closed term
of the type of natural numbers is identical to one constructed using
zero and successor only. Other examples of canonicity include the
disjunction property, which asserts that if \(P \myLogicOr Q\) is
inhabited in the empty context then either \(P\) or \(Q\) is
inhabited, and the existence property, which asserts that if
\(\myExists_{x \myElemOf A}\myApp{P}{x}\) is inhabited in the empty
context then there exists a closed term \(a \myElemOf A\) such that
\(\myApp{P}{a}\) is inhabited. Since closed terms in type theory
correspond to global sections in \((\infty,1)\)-categories, we formulate
canonicity for higher inductive types in \((\infty,1)\)-categories as
follows.

\begin{myTheorem*}[{\myRef{prop-inductive-type-canonicity}}]
  The global section functor for the initial \((\infty,1)\)-category with
  higher inductive types preserves higher inductive types.
\end{myTheorem*}

This formulation subsumes the canonicity for the type of natural numbers,
the disjunction property, and the existence property, but also yields,
for example, the following new kind of canonicity.

\begin{myTheorem*}[{\myRef{prop-homotopy-group-canonicity}}]
  Let \(C\) be the initial \((\infty,1)\)-category with higher inductive
  types. For any \(n \myGe 0\) and \(k \myGe 0\), we have a canonical
  equivalence
  \[
    \myApp{\myHoGrp_{n}}{\mySphere^{k}} \myEquiv
    \myApp{\myHom_{C}}{\myTerminal,
      \myApp{\myHoGrp_{n}}{\mySphere^{k}}}
  \]
  (of sets for \(n \myId 0\) and of groups for \(n \myGe 1\)). Here,
  the left \(\myApp{\myHoGrp_{n}}{\mySphere^{k}}\) is the \(n\)-th
  homotopy group of the \(k\)-sphere constructed in the
  \((\infty,1)\)-category of spaces, and the right
  \(\myApp{\myHoGrp_{n}}{\mySphere^{k}}\) is the one constructed in
  \(C\).
\end{myTheorem*}

We focus on \emph{finitary} higher inductive types in the sense that
all constructors are finitary operators. Colimits, the type of natural
numbers, and truncations are all finitary higher inductive
types. There are two missing important examples of higher inductive
type: W-types \myCite{martin-lof1982constructive} and localizations
\myCite{rijke2020modalities}. They are (internally) infinitary higher
inductive types. An infinitary higher inductive type \(X\) may have a
constructor that takes a function \(A \myMorphism X\) from a possibly
(internally) infinite type \(A\) as its argument. It is not hard to
extend our definition of higher inductive types to include infinitary
ones, but canonicity for infinitary higher inductive types cannot
even be stated in the same form as
\myRef{prop-inductive-type-canonicity}, because preservation of
infinitary higher inductive types for the global section functor does
not make sense since the global section functor does not preserve
function types. Because of this difficulty, it would be good to start
from finitary higher inductive types and leave it as future work to
include infinitary higher inductive types.

\myRef{sec-preliminaries} is a preliminary section. In
\myRef{sec-limit-theories} we consider higher inductive types in the
\((\infty,1)\)-category of spaces. Our take is that higher inductive types
are initial algebras for limit theories. In
\myRef{sec-inductive-types} we introduce higher inductive types in
\((\infty,1)\)-categories with finite limits and prove the finitary
presentability of the \((\infty,1)\)-category of
\((\infty,1)\)-categories with higher inductive types. We give examples of
higher inductive types in
\myRef{sec-examples-of-higher-inductive-types}. In
\myRef{sec-canonicity} we prove canonicity for higher inductive types.

\section{Preliminaries}
\label{sec-preliminaries}

We work in the language of higher category theory. By category we mean
\((\infty \myComma 1)\)-category. Every construction is homotopy
invariant. Ultimately all the results should be formalized in homotopy
type theory / univalent foundations
\myCite{hottbook,rijke2022introduction-arxiv}. With this in mind, we
use type-theoretic terminology and notation. Spaces or homotopy types
are called types. We avoid using classical axioms such as the law of
excluded middle, the axiom of choice, and Whitehead's principle
\myCite[Section 8.8]{hottbook} so that all the results are valid in an
arbitrary \(\infty\)-topos.

Every statement is implicitly parameterized by as many universes
\(U \myElemOf \myEnlarge U \myElemOf \myEnlarge^{2} U \myElemOf
\myDots\) as we need. By \myDefine{small} objects we mean objects in
\(U\). When \(X\) is a category of small objects of some kind,
\(\myEnlarge^{n} X\) denotes the category of objects in the \(n\)-th
universe of the same kind.

The reader need not be familiar with higher category theory. Most part
of higher category theory can be understood by analogy with ordinary
category theory. One thing to notice is that a category is required to
satisfy \myDefine{univalence}: for two objects \(x\) and \(y\), the
type \(x \myId y\) of identifications between \(x\) and \(y\) is
(canonically) equivalent to the type \(x \myEquiv y\) of equivalences
between \(x\) and \(y\). The following is not true in ordinary
category theory.

\begin{myProposition}
  \label{prop-mono-implies-conservative}
  If a functor \(F \myElemOf C \myMorphism D\) is monic (that is, both
  the object part and the morphism part are monic), then it is
  conservative (that is, reflects equivalences).
\end{myProposition}

\begin{myProof}
  Let \(f \myElemOf x \myMorphism y\) be a morphism in \(C\) and
  suppose that \(\myApp{F}{f}\) is an equivalence. Since the object
  part of \(F\) is monic, we have
  \((x \myId y) \myEquiv (\myApp{F}{x} \myId \myApp{F}{y})\). By
  univalence,
  \((x \myEquiv y) \myEquiv (\myApp{F}{x} \myEquiv
  \myApp{F}{y})\). Thus, \(\myApp{F}{f}\) must be the image of an
  equivalence \(f'\) in \(C\), but \(f'\) must be \(f\) as the
  morphism part of \(F\) is monic. Therefore, \(f\) is an equivalence.
\end{myProof}

\(\myCell_{1}\) denotes the walking morphism. That is, the type of
functors \(\myCell_{1} \myMorphism C\) is equivalent to the type of
morphisms in \(C\). For a category \(C\) and for a category or
simplicial type \(A\), we write \(A \myPower C\) for the category of
diagrams \(A \myMorphism C\) and call it the power of \(C\) by
\(A\). We write
\(\myDom \myComma \myCod \myElemOf \myCell_{1} \myPower C \myMorphism
C\) for the domain and codomain, respectively, projections. For
functors \(F \myElemOf C \myMorphism E\) and
\(G \myElemOf D \myMorphism E\), the comma category \((F \mySlice G)\)
is the fiber product of
\(F \myBinProd G \myElemOf C \myBinProd D \myMorphism E \myBinProd E\)
and
\((\myDom \myComma \myCod) \myElemOf \myCell_{1} \myPower E
\myMorphism E \myBinProd E\). When \(F\) is the identity functor on
\(E\), we write \((E \mySlice G)\) instead of
\((\myIdFun_{E} \mySlice G)\). In particular, \((C \mySlice x)\)
denotes the slice category over \(x\) for an object \(x \myElemOf
C\). The opposite of a category \(C\) is denoted by
\(\myApp{\myOp}{C}\). The fiber product of
\(f \myElemOf x \myMorphism z\) and \(g \myElemOf y \myMorphism z\) is
denoted by \(x \myFibBinProd{f}{g} y\) or \(x \myBinProd_{z} y\) when
\(f\) and \(g\) are clear from the context. It is also denoted by
\(\myApp{f^{\myStar}}{y}\) and called the pullback of \(y\) along
\(f\). The fiber product is regarded as an operator on
\((C \mySlice z)\) while the pullback is regarded as a functor
\(f^{\myStar} \myElemOf (C \mySlice z) \myMorphism (C \mySlice
x)\). The terminal object in a category is denoted by
\(\myTerminal\). Let \(\myType\) denote the category of small types
and let \(\myCat\) denote the category of small categories. For a
category \(C\) in \(\myEnlarge^{n} \myCat\), we write
\(\myApp{\myPsh}{C}\) for the category of presheaves over \(C\) valued
in \(\myEnlarge^{n} \myType\).

Let
\begin{equation}
  \label[square]{square-mate}
  \begin{tikzcd}
    C
    \arrow[r, "H"]
    \arrow[d, "F"']
    & C'
    \arrow[d, "F'"]
    \\ D
    \arrow[r, "K"']
    & D'
  \end{tikzcd}
\end{equation}
be a commutative square of categories and suppose that \(F\) and
\(F'\) has right adjoints \(G\) and \(G'\) respectively. The
\myDefine{mate} of \myRef{square-mate} is the composite of
natural transformations
\[
  \begin{tikzcd}
    & C
    \arrow[rr, "H"]
    \arrow[dr, "F"]
    \arrow[d, Rightarrow, "\epsilon",
    start anchor={[yshift=-1ex]},
    end anchor={[yshift=1ex]}]
    & & C'
    \arrow[rr, equal]
    \arrow[dr, "F'"']
    & {}
    \arrow[d, Rightarrow, "\eta'",
    start anchor={[yshift=-1ex]},
    end anchor={[yshift=1ex]}]
    & C'
    \\ D
    \arrow[ur, "G"]
    \arrow[rr, equal]
    & {}
    & D
    \arrow[rr, "K"']
    & & D',
    \arrow[ur, "G'"']
  \end{tikzcd}
\]
where \(\epsilon\) is the counit of the adjunction \(F \myAdjRel G\) and
\(\eta'\) is the unit of the adjunction \(F' \myAdjRel G'\). The mate of
\myRef{square-mate} when \(F\) and \(F'\) have left adjoints is dually
defined.

\subsection{Presentable categories}

We review the theory of presentable categories \myCite[Chapter
5]{lurie2009higher}. We consider the following subcategories of
\(\myEnlarge \myCat\).
\begin{itemize}
\item \(\myPrRFin\) whose objects are
  \myDefine{\(\myFinitary\)-presentable categories} (called compactly
  generated categories in \myCite{lurie2009higher}) and morphisms are
  functors preserving small limits and small filtered colimits.
\item \(\myAccFin\) whose objects are accessible categories with small
  filtered colimits and morphisms are functors preserving small
  filtered colimits.
\item \(\myAcc\) whose objects are accessible categories and morphisms
  are accessible functors.
\item \(\myPrR\) whose objects are presentable categories and
  morphisms are accessible functors preserving small limits.
\end{itemize}
We have \(\myPrRFin \mySub \myAccFin \mySub \myAcc\) and
\(\myPrRFin \mySub \myPrR \mySub \myAcc\). By the adjoint functor
theorem \myCite[Corollary 5.5.2.9]{lurie2009higher}, every morphism in
\(\myPrR\) has a left adjoint (in \(\myEnlarge \myCat\)). There are a
lot of equivalent definitions of \(\myFinitary\)-presentability of a
category \(C\), including:
\begin{itemize}
\item \(C\) is equivalent to a reflective full subcategory of
  \(\myApp{\myPsh}{A}\) closed under small filtered colimits for some
  small category \(A\);
\item \(C\) has small colimits and there exists a small full
  subcategory \(A \mySub C\) consisting of compact objects such that
  every object in \(C\) is the colimit of a small filtered diagram in
  \(A\).
\end{itemize}

We collect basic facts about presentable and accessible categories.

\begin{myProposition}
  \label{prop-acc-limit-theorem}
  \label{prop-pres-limit-theorem}
  \(\myPrRFin\), \(\myAccFin\), \(\myAcc\), and \(\myPrR\) are closed
  in \(\myEnlarge \myCat\) under small limits and powers by
  \(\myCell_{1}\).
\end{myProposition}

\begin{myProof}
  For \(\myPrR\), see \myCite[Proposition 5.5.3.6 and Theorem
  5.5.3.18]{lurie2009higher}. For \(\myAcc\), see \myCite[Proposition
  5.4.4.3 and Proposition 5.4.7.3]{lurie2009higher}. Let
  \(X \mySub \myEnlarge \myCat\) be the subcategory whose objects are
  categories with small filtered colimits and morphisms are functors
  preserving small filtered colimits. \(X\) is closed under small
  limits and powers by \(\myCell_{1}\): small filtered colimits in a
  limit or power by \(\myCell_{1}\) are computed component-wise. Since
  \(\myAccFin\) is the intersection of \(\myAcc\) and \(X\), it is
  also closed under small limits and powers by \(\myCell_{1}\). By
  \myCite[Proposition 5.5.7.6]{lurie2009higher}, \(\myPrRFin\) is
  closed under small limits. To see that it is also closed under
  powers by \(\myCell_{1}\), let \(C\) be an
  \(\myFinitary\)-presentable category which is equivalent to a
  reflective full subcategory of \(\myApp{\myPsh}{A}\) closed under
  small filtered colimits for a small category \(A\). Then
  \(\myCell_{1} \myPower C\) is equivalent to a reflective full
  subcategory of
  \(\myCell_{1} \myPower \myApp{\myPsh}{A} \myEquiv
  \myApp{\myPsh}{\myApp{\myOp}{\myCell_{1}} \myBinProd A}\) closed
  under small filtered colimits and thus \(\myFinitary\)-presentable.
\end{myProof}

\begin{myProposition}[{\myCite[cf.][Proposition 1.59]{adamek1994locally}}]
  \label{prop-lim-filt-colim-pr}
  Finite limits commute with small filtered colimits in any
  \(\myFinitary\)-presentable category. Equivalently, the limit
  functor \(\myLim \myElemOf A \myPower C \myMorphism C\) is in
  \(\myPrRFin\) for any finite simplicial type \(A\) and
  \(\myFinitary\)-presentable category \(C\).
\end{myProposition}

\begin{myProof}
  Let \(C\) be an \(\myFinitary\)-presentable category which is
  equivalent to a reflective full subcategory of \(\myApp{\myPsh}{A}\)
  closed under small filtered colimits for a small category
  \(A\). Then finite limits and small filtered colimits in \(C\) are
  computed in \(\myApp{\myPsh}{A}\). Since limits and colimits in
  \(\myApp{\myPsh}{A}\) is pointwise, the claim follows from the fact
  that finite limits and small filtered colimits commute in
  \(\myType\) \myCite[Proposition 5.3.3.3]{lurie2009higher}.
\end{myProof}

\begin{myProposition}[{\myCite[cf.][Corollary 2.47]{adamek1994locally}}]
  \label{prop-acc-lim-pres}
  An accessible category \(C\) is presentable whenever it has small
  limits.
\end{myProposition}

\begin{myProof}
  Recall that a category is presentable if it is accessible and has
  small colimits. We show that \(C\) has small colimits using the
  general adjoint functor theorem \myCite{nguyen2019adjoint}. Theorem
  4.1.1 of \myCite{nguyen2019adjoint} asserts that a functor from a
  presentable category to a category satisfying certain conditions
  that are satisfied by any accessible category is a right adjoint
  whenever it is accessible and preserves small limits. The proof of
  it uses accessibility of and small limits in the domain but no other
  consequences of presentability. It can thus be proved that a functor
  between accessible categories with small limits is a right adjoint
  whenever it is accessible and preserves small limits. Applying this
  to the diagonal functor \(C \myMorphism A \myPower C\) for a small
  simplicial type \(A\), we see that \(C\) has small colimits.
\end{myProof}

\begin{myProposition}
  \label{prop-fin-pres-closed-under-csv}
  Let \(C\) be a category with small colimits, let \(D\) be an
  \(\myFinitary\)-presentable category, and let
  \(F \myElemOf C \myMorphism D\) be a right adjoint functor
  preserving small filtered colimits. If \(F\) is conservative, then
  \(C\) is \(\myFinitary\)-presentable.
\end{myProposition}

\begin{myProof}
  \(D\) is equivalent to a reflective full subcategory of
  \(\myApp{\myPsh}{A}\) closed under small filtered colimits for some
  small category \(A\). It is thus enough to show the case when
  \(D \myId \myApp{\myPsh}{A}\). We construct a small full subcategory
  \(B \mySub C\) consisting of compact objects such that every object
  in \(C\) is the filtered colimit of a small filtered diagram in
  \(B\). Let \(G \myElemOf A \myMorphism C\) be the composite of the
  Yoneda embedding \(A \myMorphism \myApp{\myPsh}{A}\) and the left
  adjoint of \(F\). Then
  \(\myApp{\myApp{F}{y}}{x} \myEquiv \myApp{\myHom}{\myApp{G}{x}
    \myComma y}\) for \(x \myElemOf A\) and \(y \myElemOf C\). Since
  \(F\) preserves filtered colimits, every \(\myApp{G}{x}\) is
  compact. Let \(B\) be the closure of the image of \(G\) under finite
  colimits. \(B\) consists of compact objects. We show that the
  canonical morphism
  \((\myEnd^{b \myElemOf B} \myApp{\myHom}{b \myComma y} \myAct b)
  \myMorphism y\) is an equivalence for every \(y \myElemOf C\) so
  that \(y\) is the filtered colimit of the canonical diagram
  \((B \mySlice y) \myMorphism B\). Since \(F\) is conservative, it
  suffices to show that this canonical morphism is sent by \(F\) to an
  equivalence. Since \(F\) preserves small filtered colimits, this is
  equivalent to that
  \(\myEnd^{b \myElemOf B} \myApp{\myHom}{b \myComma y} \myAct
  \myApp{F}{b} \myEquiv \myApp{F}{y}\). That is,
  \(\myEnd^{b \myElemOf B} \myApp{\myHom}{b \myComma y} \myBinProd
  \myApp{\myHom}{\myApp{G}{x} \myComma b} \myEquiv
  \myApp{\myHom}{\myApp{G}{x} \myComma y}\) for every
  \(x \myElemOf A\), but this is true by Yoneda.
\end{myProof}

\subsection{Presentable \(2\)-categories}

We introduce presentable \(2\)-categories. Our definition of
\(2\)-presentability is close to Bourke's class \(\mathbf{LP}\) of
strict \(2\)-categories \myCite{bourke2021accessible} in that both are
defined in terms of accessibility of the underlying category and
\(2\)-limits. Comparison with other possible definitions in terms of
\(2\)-colimits
\myCite[cf.][]{kelly1982structures,diliberti22biaccessible-arxiv} is
beyond the scope of this paper.

Let \(\myCCat\) denote the category of small \(2\)-categories. One can
take any theory of \((\infty,2)\)-categories in the sense of Barwick and
Schommer-Pries \myCite{barwick2021unicity}. \(\myCCat\) is presentable
\myCite[Basic Data]{barwick2021unicity} and cartesian closed
\myCite[Axiom C.3]{barwick2021unicity}. We write \(\myPower\) for the
exponential in \(\myCCat\). Any category is regarded as a
\(2\)-category with trivial \(2\)-morphisms. In particular,
\(\myCCat\) is equipped with the power functor
\(C \myMapsTo \myCell_{1} \myPower C \myElemOf \myCCat \myMorphism
\myCCat\), and \(2\)-limits such as comma objects are constructed by
limits and powers by \(\myCell_{1}\). The underlying category of a
\(2\)-category \(C\) is denoted by \(\myApp{\myCoreI}{C}\).

\begin{myDefinition}
  A \(2\)-category \(C\) has \myDefine{powers by \(\myCell_{1}\)} if
  it is equipped with a \(2\)-functor \(P \myElemOf C \myMorphism C\)
  and a \(2\)-functor
  \(\epsilon \myElemOf C \myMorphism \myCell_{1} \myPower (P \mySlice C)\)
  whose projection to \(\myCell_{1} \myPower C^{2}\) is the composite
  of diagonals
  \(C \myMorphism C^{2} \myMorphism \myCell_{1} \myPower C^{2}\) such
  that, for any \(x' \myComma x \myElemOf C\), the functor
  \[
    f \myMapsTo (a \myMapsTo \myApp{\myApp{\epsilon}{x}}{a} \myComp f)
    \myElemOf \myApp{\myHHom}{x' \myComma \myApp{P}{x}} \myMorphism
    \myCell_{1} \myPower \myApp{\myHHom}{x' \myComma x}
  \]
  is an equivalence. The \(2\)-functor \(P\) is denoted by
  \(x \myMapsTo \myCell_{1} \myPower x\). A \(2\)-functor
  \(F \myElemOf C \myMorphism D\) between \(2\)-categories with powers
  by \(\myCell_{1}\) \myDefine{preserves powers by \(\myCell_{1}\)} if
  the canonical morphism
  \(\myApp{F}{\myCell_{1} \myPower x} \myMorphism \myCell_{1} \myPower
  \myApp{F}{x}\) is an equivalence for every \(x \myElemOf C\).
\end{myDefinition}

\begin{myDefinition}
  We say a \(2\)-category \(C\) is
  \myDefine{\(\myFinitary\)-presentable} if the category
  \(\myApp{\myCoreI}{C}\) is \(\myFinitary\)-presentable, \(C\) has
  powers by \(\myCell_{1}\), and the functor
  \(x \myMapsTo \myCell_{1} \myPower x \myElemOf \myApp{\myCoreI}{C}
  \myMorphism \myApp{\myCoreI}{C}\) is in \(\myPrRFin\). We write
  \(\myPPrRFin \mySub \myEnlarge \myCCat\) for the subcategory whose
  objects are the \(\myFinitary\)-presentable \(2\)-categories and
  morphisms are those \(2\)-functors \(F\) who preserve powers by
  \(\myCell_{1}\) and whose underlying functors are in \(\myPrRFin\).
\end{myDefinition}

\begin{myProposition}
  \label{prop-2pr-limit-theorem}
  \(\myPPrRFin \mySub \myEnlarge \myCCat\) is closed under small
  limits and powers by \(\myCell_{1}\).
\end{myProposition}

\begin{myProof}
  Let \(C \myElemOf I \myMorphism \myPPrRFin\) be a small diagram and
  let \(D\) be its limit in \(\myEnlarge \myCCat\). By
  \myRef{prop-pres-limit-theorem}, \(\myApp{\myCoreI}{D}\) is in
  \(\myPrRFin\). Observe that powers by \(\myCell_{1}\) in \(D\) are
  computed component-wise. Since small filtered colimits in
  \(\myApp{\myCoreI}{D}\) are also computed component-wise, it follows
  that \(D\) is \(\myFinitary\)-presentable. \(\myPPrRFin\) is closed
  under powers by \(\myCell_{1}\) for the same reason.
\end{myProof}

\(\myCat\) is extended to a \(2\)-category such that
\(\myApp{\myHom_{\myCat}}{C \myComma D} \myEquiv C \myPower D\). We
refer the reader to \myCite{lurie2009goodwillie} for an explicit
construction of \(\myCat\).

\begin{myProposition}
  \label{prop-cat-2pr}
  \(\myCat\) is extended to an \(\myFinitary\)-presentable
  \(2\)-category.
\end{myProposition}

\begin{myProof}
  The category \(\myCat\) is \(\myFinitary\)-presentable because it is
  the reflective full subcategory of the category
  \(\myApp{\myPsh}{\mySimpCat}\) of simplicial types consisting of
  complete Segal objects and is closed under small filtered colimits. It
  remains to show that
  \(C \myMapsTo \myCell_{1} \myPower C \myElemOf \myCat \myMorphism
  \myCat\) preserves small filtered colimits. Since small filtered
  colimits in \(\myCat\) are computed in
  \(\myApp{\myPsh}{\mySimpCat}\) and since \(\myCell_{1} \myPower C\)
  is the exponential \(C^{\myStdSimp{1}}\) in
  \(\myApp{\myPsh}{\mySimpCat}\), it suffices to show that the
  exponential by \(\myStdSimp{1}\) preserves small filtered
  colimits. But this follows from the fact that the cartesian product
  with \(\myStdSimp{1}\) preserves finite simplicial types.
\end{myProof}

\begin{myProposition}
  \label{prop-limit-2-functor-2pr}
  Let \(C\) be an \(\myFinitary\)-presentable \(2\)-category. The power
  \(2\)-functor
  \(x \myMapsTo \myCell_{1} \myPower x \myElemOf C \myMorphism C\) is
  in \(\myPPrRFin\), and for any finite simplicial type \(A\), the
  limit \(2\)-functor \(\myLim \myElemOf A \myPower C \myMorphism C\)
  is in \(\myPPrRFin\).
\end{myProposition}

\begin{myProof}
  This is because limits commute with powers by \(\myCell_{1}\).
\end{myProof}

By \myRef{prop-2pr-limit-theorem}, \myRef{prop-cat-2pr}, and
\myRef{prop-limit-2-functor-2pr}, the \(2\)-category of small
instances of a structure defined by finite \(2\)-limits is
\(\myFinitary\)-presentable. Here, a structure is defined by finite
\(2\)-limits when it is defined by the following specifications.
\begin{enumerate}[series=pres-specification, label=(\alph*), ref=(\alph*)]
\item \label[specification]{spec-categorical-first} Introducing a new object.
\item Introducing a new morphism between objects constructed
  from previously introduced data using finite \(2\)-limits.
\item Introducing a new identification between morphisms constructed
  from previously introduced data using finite \(2\)-limits.
  \label[specification]{spec-categorical-last}
\end{enumerate}
Here, ``using finite \(2\)-limits'' means ``using finite limits and
powers by \(\myCell_{1}\)''. These specifications correspond to
constructions in \(\myPPrRFin\) as in \myRef{tab-spec-pres}. The
following specifications are combinations of specifications
\myLabelRef{spec-categorical-first}--\myLabelRef{spec-categorical-last}.
\begin{enumerate}[resume*=pres-specification]
\item Introducing a new \(2\)-morphism between morphisms constructed
  from previously introduced data using finite \(2\)-limits. (A
  \(2\)-morphism between morphisms
  \(f \myComma g \myElemOf x \myMorphism y\) is represented by a
  morphism \(t \myElemOf x \myMorphism \myCell_{1} \myPower y\) and
  identifications \(\myDom \myComp t \myId f\) and
  \(\myCod \myComp t \myId g\).)
\item Introducing a new identification between \(2\)-morphisms
  constructed from previously introduced data using finite
  \(2\)-limits.
\item Requiring a morphism constructed from previously introduced data
  using finite \(2\)-limits is an equivalence/ left adjoint / right
  adjoint.
\item Requiring a morphism constructed from previously introduced data
  using finite \(2\)-limits is fully faithful. (A morphism
  \(f \myElemOf x \myMorphism y\) is fully faithful if
  \(\myCell_{1} \myPower x \myMorphism (\myCell_{1} \myPower y)
  \myBinProd_{y^{2}} (x^{2})\) is an equivalence.)
\end{enumerate}
For example, we have the \(\myFinitary\)-presentable \(2\)-category
\(\myLex\) of small categories with finite limits. Indeed, finite
limits are a structure defined by finite \(2\)-limits since a finite
limit is the right adjoint of the diagonal functor
\(C \myMorphism A \myPower C\) for a finite simplicial type \(A\).

\begin{table}
  \begin{tabular}{ll}
    Specification
    & Construction
    \\ \hline
    object
    & \(D \myDefEq C \myBinProd \myCat\)
    \\ morphism \(x \myMorphism y\)
    & pullback \(
      \begin{tikzcd}
        D
        \arrow[r]
        \arrow[d]
        & \myCell_{1} \myPower \myCat
        \arrow[d, "{(\myDom \myComma \myCod)}"]
        \\ C
        \arrow[r, "{(x \myComma y)}"']
        & \myCat^{2}
      \end{tikzcd}
      \)
    \\ identification \(f \myId g\)
    & pullback \(
      \begin{tikzcd}
        D
        \arrow[r]
        \arrow[d]
        & \myCell_{1} \myPower \myCat
        \arrow[d, "\myDiagonal"]
        \\ C
        \arrow[r, "{(f \myComma g)}"']
        & (\myCell_{1} \myPower \myCat) \myBinProd_{\myCat^{2}}
        (\myCell_{1} \myPower \myCat)
      \end{tikzcd}
      \)
  \end{tabular}
  \caption{Specification of structure and construction in
    \(\myPPrRFin\)}
  Suppose we have defined some structure \(X\) and let \(C\) be the
  \(2\)-category of small instances of \(X\). When \(Y\) is an
  extension of \(X\) by one of data on the left, the category \(D\) of
  small instances of \(Y\) is obtained by the construction on the
  right in the same row.
  \label{tab-spec-pres}
\end{table}

For a structure defined by finite \(2\)-limits, we will not explicitly
define what a morphism for that structure is. Rather, the notion of
morphism is derived from the construction in \(\myPPrRFin\). For all
the structures introduced in this paper, a morphism is a
structure-preserving functor.

We end this subsection with a few facts in \(2\)-category theory.

\begin{myProposition}
  \label{prop-2pr-left-2-adjoint}
  Let \(F \myElemOf D \myMorphism C\) be a morphism in
  \(\myPPrRFin\). Let \(G\) be the left adjoint of
  \(\myApp{\myCoreI}{F}\) with unit \(\eta\). Then, for any
  \(x \myElemOf C\) and \(y \myElemOf D\), the functor
  \[
    \myApp{\myHHom}{\myApp{G}{x} \myComma y}
    \myXMorphism{F}
    \myApp{\myHHom}{\myApp{F}{\myApp{G}{x}} \myComma \myApp{F}{y}}
    \myXMorphism{\eta^{\myStar}}
    \myApp{\myHHom}{x \myComma \myApp{F}{y}}
  \]
  is an equivalence.
\end{myProposition}

\begin{myProof}
  We know that the object part of the functor is an equivalence. The
  morphism part
  \[
    \myApp{\myHom}{\myCell_{1} \myComma \eta^{\myStar} \myComp F}
    \myElemOf \myApp{\myHom}{\myCell_{1} \myComma
      \myApp{\myHHom}{\myApp{G}{x} \myComma y}} \myMorphism
    \myApp{\myHom}{\myCell_{1} \myComma \myApp{\myHHom}{x \myComma
        \myApp{F}{y}}}
  \]
  is equivalent to the equivalence
  \(\myApp{\myObj}{\myApp{\myHHom}{\myApp{G}{x} \myComma \myCell_{1}
      \myPower y}} \myEquiv \myApp{\myObj}{\myApp{\myHHom}{x \myComma
      \myApp{F}{\myCell_{1} \myPower y}}}\) because \(F\) preserves
  powers by \(\myCell_{1}\).
\end{myProof}

\begin{myConstruction}
  \label{cst-2pr-left-2-functor}
  Let \(F\) be a morphism in \(\myPPrRFin\). By
  \myRef{prop-2pr-left-2-adjoint}, the morphism part of the left
  adjoint \(G\) of \(\myApp{\myCoreI}{F}\) is extended to a functor as
  \(\myApp{\myHHom}{x \myComma x'} \myXMorphism{\eta_{\myStar}}
  \myApp{\myHHom}{x \myComma \myApp{F}{\myApp{G}{x'}}} \myEquiv
  \myApp{\myHHom}{\myApp{G}{x} \myComma \myApp{G}{x'}}\).
\end{myConstruction}

\begin{myDefinition}
  Let \(C\) be a \(2\)-category. We say a sub-\(2\)-category
  \(C' \mySub C\) is \myDefine{locally full} if
  \(\myApp{\myHHom_{C'}}{x \myComma x'} \mySub \myApp{\myHHom_{C}}{x
    \myComma x'}\) is full for all \(x \myComma x' \myElemOf C'\). A
  locally full sub-\(2\)-category \(C' \mySub C\) is \myDefine{locally
    closed under retract} if
  \(\myApp{\myHHom_{C'}}{x \myComma x'} \mySub \myApp{\myHHom_{C}}{x
    \myComma x'}\) is closed under retract for all
  \(x \myComma x' \myElemOf C'\).
\end{myDefinition}

\begin{myProposition}
  \label{prop-adj-locally-closed-under-retract}
  Let \(\myAdj\) denote the \(\myFinitary\)-presentable \(2\)-category
  of adjunctions. Both the left and right projections
  \(\myAdj \myMorphism \myCell_{1} \myPower \myCat\) are mono and
  exhibit \(\myAdj\) as locally full sub-\(2\)-categories of
  \(\myCell_{1} \myPower \myCat\) which are locally closed under
  retract.
\end{myProposition}

\begin{myProof}
  We prove only for the left projection. The right projection is
  dual. The left projection
  \(\myAdj \myMorphism \myCell_{1} \myPower \myCat\) induces an
  equivalence \(\myAdj \myEquiv \myLAdj\), where \(\myLAdj\) is the
  locally full sub-\(2\)-category of \(\myCell_{1} \myPower \myCat\)
  whose objects are those functors who have right adjoints and
  morphisms are those squares whose mates are equivalences. Let
  \((H \myComma K) \myElemOf F \myMorphism F'\) be a morphism in
  \(\myLAdj\) and let \((H' \myComma K')\) be a retract of
  \((H \myComma K)\) in
  \(\myApp{\myHHom_{\myCell_{1} \myPower \myCat}}{F \myComma
    F'}\). One can see that the mate of \((H' \myComma K')\) is a
  retract of the mate of \((H \myComma K)\) and thus an equivalence.
\end{myProof}

\subsection{Lex universes}

We review univalent universes in categories with finite limits. See
\myCite{gepner2017univalence,rasekh2021univalence,nguyen2022type-arxiv}
for more details.

Let \(C\) be a category with finite limits. A morphism
\(u \myElemOf y \myMorphism x\) in \(C\) is \myDefine{exponentiable}
if the following equivalent conditions hold:
\begin{enumerate}
\item \label{cond-exponentiable-pushforward} The pullback functor
  \(u^{\myStar} \myElemOf (C \mySlice x) \myMorphism (C \mySlice y)\)
  has a right adjoint \(u_{\myStar}\).
\item \label{cond-exponentiable-local-exp} The fiber product functor
  \(z \myMapsTo y \myBinProd_{x} z \myElemOf (C \mySlice x)
  \myMorphism (C \mySlice x)\) has a right adjoint
  \(z \myMapsTo y \myExp_{x} z\).
\item \label{cond-exponentiable-fam} The fiber product functor
  \(z \myMapsTo y \myBinProd_{x} z \myElemOf (C \mySlice x)
  \myMorphism C\) has a right adjoint \(\myFam_{u}\).
\end{enumerate}
Moreover, \(u_{\myStar}\), \(z \myMapsTo y \myExp_{x} z\), and
\(\myFam_{u}\) are constructible from each other:
\(y \myExp_{x} z \myEquiv
\myApp{u_{\myStar}}{\myApp{u^{\myStar}}{z}}\);
\(\myApp{\myFam_{u}}{z} \myEquiv y \myExp_{x} (z \myBinProd x)\);
\(\myApp{u_{\myStar}}{z} \myEquiv
\myApp{\eta_{x}^{\myStar}}{\myApp{\myFam_{u}}{z}}\) where
\(\eta_{x} \myElemOf x \myMorphism \myApp{\myFam_{u}}{y}\) is the unit
of \(\myFam_{u}\). Exponentiable morphisms are closed under pullback:
if \(u' \myElemOf y' \myMorphism x'\) is the pullback of \(u\) along a
morphism \(f \myElemOf x' \myMorphism x\), then
\(\myApp{\myFam_{u'}}{z} \myEquiv
\myApp{f^{\myStar}}{\myApp{\myFam_{u}}{z}}\). For two exponentiable
morphisms \(u \myElemOf y \myMorphism x\) and
\(u' \myElemOf y' \myMorphism x'\), we define a morphism
\(u \myPolyProd u'\) in \(C\) to be the composite
\(\myApp{\epsilon^{\myStar}}{y'} \myMorphism \myApp{\myFam_{u}}{x'}
\myBinProd_{x} y \myXMorphism{\myProj_{1}} \myApp{\myFam_{u}}{x'}\),
where
\(\epsilon \myElemOf \myApp{\myFam_{u}}{x'} \myBinProd_{x} y \myMorphism x'\)
is the counit of the right adjoint \(\myFam_{u}\). Since
\(u \myPolyProd u'\) is the composite of pullbacks of exponentiable
morphisms, it is exponentiable.

Let \(u \myElemOf y \myMorphism x\) be an exponentiable morphism in
\(C\). We define
\(\myApp{\myIArr}{u} \myElemOf (C \mySlice x \myBinProd x)\) to be
\(\myApp{\myProj_{1}^{\myStar}}{y} \myExp_{x \myBinProd x}
\myApp{\myProj_{2}^{\myStar}}{y}\). For
\((f \myComma g) \myElemOf z \myMorphism x \myBinProd x\), a morphism
\(z \myMorphism \myApp{\myIArr}{u}\) over \(x \myBinProd x\)
corresponds to a morphism
\(\myApp{f^{\myStar}}{y} \myMorphism \myApp{g^{\myStar}}{y}\) over
\(z\). Let \(\myApp{\myIArrTwo}{u}\) be the fiber product of
\(\myProj_{2} \myElemOf \myApp{\myIArr}{u} \myMorphism x\) and
\(\myProj_{1} \myElemOf \myApp{\myIArr}{u} \myMorphism x\). We have
the composition morphism
\(\myIComp \myElemOf \myApp{\myIArrTwo}{u} \myMorphism
\myApp{\myIArr}{u}\) and the identity morphism
\(\myIIdFun \myElemOf x \myMorphism \myApp{\myIArr}{u}\). Let
\(\myApp{\myIRet}{u}\) be the fiber product of \(\myIComp\) and
\(\myIIdFun\). We define
\(\myApp{\myIArrEq}{u} \myElemOf (C \mySlice \myApp{\myIArr}{u})\) to
be the fiber product of
\(\myApp{\myIRet}{u} \myMorphism \myApp{\myIArrTwo}{u}
\myXMorphism{\myProj_{1}} \myApp{\myIArr}{u}\) and
\(\myApp{\myIRet}{u} \myMorphism \myApp{\myIArrTwo}{u}
\myXMorphism{\myProj_{2}} \myApp{\myIArr}{u}\). The morphism
\(\myApp{\myIArrEq}{u} \myMorphism \myApp{\myIArr}{u}\) is monic, and a
morphism \(z \myMorphism \myApp{\myIArr}{u}\) over \(x \myBinProd x\)
factors through \(\myApp{\myIArrEq}{u}\) if and only if the
corresponding morphism
\(\myApp{f^{\myStar}}{y} \myMorphism \myApp{g^{\myStar}}{y}\) over
\(z\) is an equivalence. In particular,
\(\myIIdFun \myElemOf x \myMorphism \myApp{\myIArr}{u}\) factors
through \(\myApp{\myIArrEq}{u}\). We say \(u\) is \myDefine{univalent}
if the morphism
\(\myIIdFun \myElemOf x \myMorphism \myApp{\myIArrEq}{u}\) is an
equivalence.

When \(u \myElemOf y \myMorphism x\) is univalent, any morphism in
\(C\) is a pullback of \(u\) in at most one way. We say a morphism
\(u' \myElemOf y' \myMorphism x'\) in \(C\) is \myDefine{\(u\)-small}
if it is a pullback of \(u\). If this is the case, the morphism
\(x' \myMorphism x\) is called the characteristic morphism of
\(u'\). We say:
\begin{itemize}
\item \(u\) \myDefine{supports a unit type} if the identity on
  \(\myTerminal\) is \(u\)-small;
\item \(u\) \myDefine{supports pair types} if \(u \myPolyProd u\) is
  \(u\)-small;
\item \(u\) \myDefine{supports identity types} if the diagonal
  morphism \(y \myMorphism y \myBinProd_{x} y\) is \(u\)-small.
\end{itemize}

\begin{myDefinition}
  By a \myDefine{lex universe} we mean a univalent exponentiable
  morphism that supports a unit type, pair types, and identity
  types. Dually, a \myDefine{co-lex-universe} is defined as a lex
  universe in the opposite category. For a lex universe
  \(u \myElemOf y \myMorphism x\), we refer to the characteristic
  morphisms of the identity on \(\myTerminal\) and the diagonal
  \(y \myMorphism y \myBinProd_{x} y\) as
  \(\myIntern{\myTerminal} \myElemOf \myTerminal \myMorphism x\) and
  \(\myIntern{\myId} \myElemOf y \myBinProd_{x} y \myMorphism x\),
  respectively.
\end{myDefinition}

\begin{myProposition}
  \label{prop-small-morphisms-cancel}
  Let \(C\) be a category with finite limits, let
  \(u \myElemOf y \myMorphism x\) be a lex universe in \(C\), and let
  \(u' \myElemOf y' \myMorphism x'\) and
  \(f \myElemOf y'' \myMorphism y'\) be morphisms in \(C\). Suppose
  that \(u'\) is \(u\)-small. Then \(f\) is \(u\)-small if and only if
  \(u' \myComp f\) is.
\end{myProposition}

\begin{myProof}
  The ``only if'' part is true because \(u\) supports pair types. The
  ``if'' part follows from \myCite[Lemma 5.19]{nguyen2022type-arxiv}.
\end{myProof}

\subsection{Marked categories}

\begin{myDefinition}
  A \myDefine{marked category} is a category \(C\) equipped with a
  full subcategory
  \(\myCell_{1} \myMkPower C \mySub \myCell_{1} \myPower C\)
  containing all the identities. Objects in
  \(\myCell_{1} \myMkPower C\) are called \myDefine{marked} morphisms
  in \(C\). By definition, we have the \(\myFinitary\)-presentable
  \(2\)-category \(\myMkCat\) of small marked categories. Morphisms in
  \(\myMkCat\) are called \myDefine{marked functors}.
\end{myDefinition}

\begin{myDefinition}
  A marked functor \(P \myElemOf T \myMorphism C\) between marked
  categories is a \myDefine{marked right fibration} if the functor
  \[
    (\myCell_{1} \myPower P \myComma \myCod) \myElemOf \myCell_{1}
    \myMkPower T \myMorphism (\myCell_{1}
    \myMkPower C) \myFibBinProd{\myCod}{P} T
  \]
  is an equivalence. We have the \(\myFinitary\)-presentable
  \(2\)-category \(\myMkRFib\) of small marked right fibrations.
\end{myDefinition}

The following are immediate from the definition.

\begin{myProposition}
  \label{prop-marked-rfib-cancel}
  Let \(F \myElemOf S \myMorphism T\) and
  \(P \myElemOf T \myMorphism C\) be marked functors between marked
  categories. Suppose that \(P\) is a marked right fibration. Then
  \(F\) is a marked right fibration if and only if \(P \myComp F\)
  is. \myQED
\end{myProposition}

\begin{myProposition}
  \label{prop-marked-rfib-pullback}
  Marked right fibrations are closed under pullback along an arbitrary
  marked functor. \myQED
\end{myProposition}

A useful property of marked right fibrations is
\myRef{prop-lift-marked-adjoint} below which asserts that, under some
mild assumptions, a right adjoint lifts along a marked right
fibration.

\begin{myDefinition}
  Let \(P \myElemOf T \myMorphism C\) be a functor. We say a morphism
  \(f \myElemOf t' \myMorphism t\) is \myDefine{\(P\)-cartesian} if
  for any object \(t'' \myElemOf T\), the square
  \[
    \begin{tikzcd}
      \myApp{\myHom}{t'' \myComma t'}
      \arrow[r, "f_{\myStar}"]
      \arrow[d, "P"']
      & \myApp{\myHom}{t'' \myComma t}
      \arrow[d, "P"]
      \\ \myApp{\myHom}{\myApp{P}{t''} \myComma \myApp{P}{t'}}
      \arrow[r, "\myApp{P}{f}_{\myStar}"']
      & \myApp{\myHom}{\myApp{P}{t''} \myComma \myApp{P}{t}}
    \end{tikzcd}
  \]
  is a pullback.
\end{myDefinition}

\begin{myLemma}
  \label{prop-marked-is-cartesian}
  Let \(P \myElemOf T \myMorphism C\) be a marked right fibration
  between marked categories. Then a morphism \(f\) in \(T\) is marked
  if and only if \(\myApp{P}{f}\) is marked and \(f\) is
  \(P\)-cartesian.
\end{myLemma}

\begin{myProof}
  We first show that every marked morphism
  \(f \myElemOf t' \myMorphism t\) in \(T\) is \(P\)-cartesian. Let
  \(t'' \myElemOf T\). Because \(\myDom\) is the right adjoint of the
  diagonal functor \(C \myMorphism \myCell_{1} \myPower C\), we have
  \(\myApp{\myHom}{t'' \myComma t'} \myEquiv
  \myApp{\myHom_{\myCell_{1} \myMkPower T}}{\myIdFun_{t''} \myComma
    f}\). Consider the following diagram.
  \[
    \begin{tikzcd}
      \myApp{\myHom}{t'' \myComma t'}
      \arrow[r]
      \arrow[d, "\myEquiv"']
      & \myApp{\myHom}{\myApp{P}{t''} \myComma \myApp{P}{t'}}
      \myBinProd_{\myApp{\myHom}{\myApp{P}{t''} \myComma
          \myApp{P}{t}}} \myApp{\myHom}{t'' \myComma t}
      \arrow[d, "\myEquiv"]
      \\ \myApp{\myHom_{\myCell_{1} \myMkPower T}}{\myIdFun_{t''}
        \myComma f}
      \arrow[r]
      & \myApp{\myHom_{\myCell_{1} \myMkPower
          C}}{\myIdFun_{\myApp{P}{t''}} \myComma \myApp{P}{f}}
      \myBinProd_{\myApp{\myHom}{\myApp{P}{t''} \myComma
          \myApp{P}{t}}} \myApp{\myHom}{t'' \myComma t}
    \end{tikzcd}
  \]
  The bottom map is an equivalence since \(P\) is a marked right
  fibration. Hence, the top map is also an equivalence, which means
  that \(f\) is \(P\)-cartesian.

  Conversely, suppose that \(\myApp{P}{f}\) is marked and that \(f\)
  is \(P\)-cartesian. Since \(P\) is a marked right fibration, there
  exists a unique pair \((s \myComma g)\) of \(s \myElemOf T\) and
  marked \(g \myElemOf s \myMorphism t\) sent by \(P\) to
  \((\myApp{P}{t'} \myComma \myApp{P}{f})\). As we have proved, \(g\)
  is \(P\)-cartesian. By the definition of \(P\)-cartesianness, we
  have a natural equivalence
  \(\myApp{\myHom}{t'' \myComma t'} \myEquiv \myApp{\myHom}{t''
    \myComma s}\) that commutes with \(f_{\myStar}\) and
  \(g_{\myStar}\) for all \(t'' \myElemOf T\). By Yoneda, it follows
  that \((t' \myComma f) \myEquiv (s \myComma g)\). Since \(g\) is
  marked, so is \(f\).
\end{myProof}

\begin{myProposition}
  \label{prop-lift-marked-adjoint}
  Let
  \[
    \begin{tikzcd}
      T'
      \arrow[r, "\bar{F}"]
      \arrow[d, "P'"']
      & T
      \arrow[d, "P"]
      \\ C'
      \arrow[r, "F"']
      & C
    \end{tikzcd}
  \]
  be a pullback of categories where \(T\) and \(C\) are marked
  categories and \(P\) is a marked right fibration. Suppose that \(F\)
  has a right adjoint \(G\) whose counit
  \(\epsilon \myElemOf F \myComp G \myMMorphism \myIdFun\) is component-wise
  marked. Then \(\bar{F}\) has a right adjoint \(\bar{G}\) whose
  counit \(\bar{\epsilon}\) is component-wise marked, and the mate
  \(P' \myComp \bar{G} \myMMorphism G \myComp P\) is an equivalence.
\end{myProposition}

\begin{myProof}
  Let \(t \myElemOf T\) and let \(x \myDefEq \myApp{P}{t}\). Since
  \(P\) is a marked right fibration and since
  \(\epsilon_{x} \myElemOf \myApp{F}{\myApp{G}{x}} \myMorphism x\) is marked,
  there exists a unique pair
  \((\myApp{\bar{G}}{t} \myComma \bar{\epsilon}_{t})\) of
  \(\myApp{\bar{G}}{t} \myElemOf T\) and marked
  \(\bar{\epsilon}_{t} \myElemOf \myApp{\bar{G}}{t} \myMorphism t\) sent by
  \(P\) to \((\myApp{F}{\myApp{G}{x}} \myComma \epsilon_{x})\). The pair
  \((\myApp{G}{x} \myComma \myApp{\bar{G}}{t})\) defines an object in
  \(T'\) and \(\bar{\epsilon}_{t}\) is a morphism
  \(\myApp{\bar{F}}{\myApp{G}{x} \myComma \myApp{\bar{G}}{t}}
  \myMorphism t\) in \(T\). For \((x' \myComma t') \myElemOf T'\),
  consider the following diagram.
  \[
    \begin{tikzcd}
      \myApp{\myHom}{(x' \myComma t') \myComma (\myApp{G}{x}
        \myComma \myApp{\bar{G}}{t})}
      \arrow[r, "\bar{F}"]
      \arrow[d, "P'"']
      & \myApp{\myHom}{t' \myComma \myApp{\bar{G}}{t}}
      \arrow[r, "\bar{\epsilon}_{\myStar}"]
      \arrow[d, "P"]
      & \myApp{\myHom}{t' \myComma t}
      \arrow[d, "P"]
      \\ \myApp{\myHom}{x' \myComma \myApp{G}{x}}
      \arrow[r, "F"']
      & \myApp{\myHom}{\myApp{F}{x'} \myComma \myApp{F}{\myApp{G}{x}}}
      \arrow[r, "\epsilon_{\myStar}"']
      & \myApp{\myHom}{\myApp{F}{x'} \myComma x}
    \end{tikzcd}
  \]
  The left square is a pullback by definition. The right square is a
  pullback by \myRef{prop-marked-is-cartesian}. The composite of the
  bottom maps is an equivalence as \(F \myAdjRel G\). Therefore, the
  composite of the top maps is also an equivalence, and thus
  \(\bar{G}\) and \(\bar{\epsilon}\) define a right adjoint of
  \(\bar{F}\). By construction, \(\bar{\epsilon}\) is component-wise marked
  and \(P' \myComp \bar{G} \myEquiv G \myComp P\).
\end{myProof}

\begin{myConstruction}
  \label{cst-marked-rfib-2-functorial}
  Let \(S \myMorphism T\) be a marked right fibration, let
  \(F \myComma G \myElemOf T' \myMorphism T\) be marked functors, and
  let \(p \myElemOf F \myMMorphism G\) be a natural transformation
  that is component-wise marked. We define a marked functor
  \(\myApp{p^{\myStar}}{S} \myElemOf \myApp{G^{\myStar}}{S}
  \myMorphism \myApp{F^{\myStar}}{S}\) by sending
  \((t' \myComma s) \myElemOf \myApp{G^{\myStar}}{S}\) to
  \((t' \myComma \myApp{p_{t'}^{\myStar}}{s})\) where
  \(\myApp{p_{t'}^{\myStar}}{s} \myMorphism s\) is the unique marked
  morphism over
  \(p_{t'} \myElemOf \myApp{F}{t'} \myMorphism \myApp{G}{t'}\).
\end{myConstruction}

\subsection{Cartesian fibrations}

\begin{myConstruction}
  For any category \(C\), we have a marked category
  \(\myMaxMark C\) whose underlying category is \(C\) and marked
  morphisms are all the morphisms. This defines a morphism
  \(\myMaxMark \myElemOf \myCat \myMorphism \myMkCat\) in
  \(\myPPrRFin\). We define the \(\myFinitary\)-presentable
  \(2\)-category \(\myCart\) of \myDefine{cartesian fibrations} to be
  the fiber product of
  \(\myMaxMark \myElemOf \myCat \myMorphism \myMkCat\) and
  \(\myCod \myElemOf \myMkRFib \myMorphism \myMkCat\).
\end{myConstruction}

Unwinding the definition, a cartesian fibration consists of a category
\(C\), a marked category \(T\), and a marked right fibration
\(P \myElemOf T \myMorphism \myMaxMark C\). By
\myRef{prop-marked-is-cartesian}, the class of marked morphisms in
\(T\) must be the class of the \(P\)-cartesian morphisms. Thus, our
definition of cartesian fibrations coincides with the usual definition
by the existence of cartesian lifts, and \(\myCart\) is regarded as a
subcategory of \(\myCell_{1} \myPower \myCat\). We also regard
\(\myCart\) as a cartesian fibration over \(\myCat\) by the codomain
projection \(\myCod \myElemOf \myCart \myMorphism \myCat\) where the
cartesian lift is given by pullback.

The straightening-unstraightening theorem \myCite[Theorem
3.2.0.1]{lurie2009higher} asserts that cartesian fibrations over \(C\)
with small fibers are equivalent to functors
\(\myApp{\myOp}{C} \myMorphism \myCat\). The left-to-right
construction maps a cartesian fibration \(T\) over \(C\) to the
functor that assigns the fiber \(\myFiber{T}{x}\) for every
\(x \myElemOf C\). Many category-theoretic operations then induce
operations on cartesian fibrations over \(C\). We use subscript
\({}_{C}\) to indicate operations on cartesian fibrations over
\(C\). For example, \(\myApp{\myOp_{C}}{T}\) is the cartesian
fibration over \(C\) such that
\(\myFiber{\myApp{\myOp_{C}}{T}}{x} \myEquiv
\myApp{\myOp}{\myFiber{T}{x}}\). Another example is
\(A \myPower_{C} T\) for a simplicial type \(A\) which is defined by
\(\myFiber{(A \myPower_{C} T)}{x} \myEquiv A \myPower \myFiber{T}{x}\)
or as the pullback of \(A \myPower T\) along the diagonal functor
\(C \myMorphism A \myPower C\). The functor
\((C \myComma T) \myMapsTo (C \myComma A \myPower_{C} T) \myElemOf
\myCart \myMorphism \myCart\) is in \(\myPPrRFin\) whenever \(A\) is
finite. The comma \((F \mySlice_{C} G)\) for
\(F \myElemOf T' \myMorphism T\) and \(G \myElemOf T'' \myMorphism T\)
over \(C\) is defined to be
\((T' \myBinProd_{C} T'') \myBinProd_{T \myBinProd_{C} T} (\myCell_{1}
\myPower_{C} T)\).

A structure \(X\) defined by finite \(2\)-limits is translated into a
structure \(X_{\myCartMark}\) defined by finite \(2\)-limits as
follows.
\begin{enumerate}
\item Introduce an object \({*}\).
\item For every object \(x\) introduced in \(X\), introduce a
  cartesian fibration \(x'\) over \({*}\).
\item For every morphism \(f \myElemOf x \myMorphism y\) introduced in
  \(X\), introduce a functor \(f' \myElemOf x' \myMorphism y'\) over
  \({*}\) preserving cartesian morphisms, where \(x'\) and \(y'\) are
  obtained from the constructions of \(x\) and \(y\), respectively, by
  replacing finite limits by finite limits over \({*}\) and
  \(\myCell_{1} \myPower\) by \(\myCell_{1} \myPower_{{*}}\).
\item For every identification \(p \myElemOf f \myId g\) introduced in
  \(X\), introduce an identification \(p' \myElemOf f' \myId g'\) over
  \({*}\), where \(f'\)and \(g'\) are obtained from the constructions
  of \(f\) and \(g\), respectively, by replacing finite limits by
  finite limits over \({*}\) and \(\myCell_{1} \myPower\) by
  \(\myCell_{1} \myPower_{{*}}\).
\end{enumerate}
The structure \(X_{\myCartMark}\) is called \myDefine{cartesian
  \(X\)}. For example: a \myDefine{cartesian functor} consists of a
category \(C\), cartesian fibrations \(T\) and \(T'\) over \(C\), and
a functor \(F \myElemOf T \myMorphism T'\) over \(C\) preserving
cartesian morphisms; a \myDefine{cartesian adjunction} consists of a
category \(C\), cartesian fibrations \(T\) and \(T'\), and an
adjunction \(F \myAdjRel G\) between \(T\) and \(T'\) over \(C\) such
that both \(F\) and \(G\) are cartesian functors; a cartesian
fibration \(T\) over \(C\) has \myDefine{finite cartesian colimits} if
the diagonal functor \(T \myMorphism A \myPower_{C} T\) has a
cartesian left adjoint over \(C\) for every finite simplicial type
\(A\).

\begin{myDefinition}
  Let \(S\) and \(T\) be cartesian fibrations over \(C\) and let
  \(F \myElemOf S \myMorphism T\) be a cartesian functor over
  \(C\). We say \(F\) is \myDefine{representable} if it reflects
  cartesian morphisms and has a (non-cartesian) right adjoint
  \(F^{\myRepRA}\) whose counit is component-wise cartesian.
\end{myDefinition}

When \(F\) is representable, the unit is also component-wise cartesian
and \(F^{\myRepRA}\) sends cartesian morphisms to cartesian morphisms
(so \(F\) is a generalized category with families in the sense of
Coraglia and Emmenegger \myCite{coraglia2024analysis}). By definition,
representable cartesian functors are closed under
composition. Representable cartesian functors are closed under
pullback in the sense that if
\[
  \begin{tikzcd}
    S'
    \arrow[r, "G"]
    \arrow[d, "F'"']
    & S
    \arrow[d, "F"]
    \\ T'
    \arrow[r, "H"']
    & T
  \end{tikzcd}
\]
is a pullback in \(\myFiber{\myCart}{C}\) and \(F\) is representable,
then \(F'\) is representable and the mate
\(G \myComp (F')^{\myRepRA} \myMMorphism F^{\myRepRA} \myComp H\) is
an equivalence. This follows from
\myRef{prop-lift-marked-adjoint}.

\begin{myLemma}
  \label{prop-cart-fib-pullback-by-ext}
  Let \(F \myElemOf S \myMorphism T\) be a cartesian functor between
  cartesian fibrations over \(C\), let \(T'\) be a cartesian fibration
  over \(C'\), and let
  \[
    \begin{tikzcd}
      S'
      \arrow[r, "G"]
      \arrow[d, "F'"']
      & S
      \arrow[d, "F"]
      \\ T'
      \arrow[r, "H"']
      & T
    \end{tikzcd}
  \]
  be a pullback of categories. Suppose that \(H\) sends cartesian
  morphisms to cartesian morphisms. Then the composite \(S'
  \myXMorphism{F'} T' \myMorphism C'\) is a cartesian fibration and
  \(F'\) is a cartesian functor.
\end{myLemma}

\begin{myProof}
  We regard the given pullback square as a pullback in \(\myMkCat\)
  where the marked morphisms in \(S\), \(T\), and \(T'\) are the
  cartesian morphisms. By \myRef{prop-marked-rfib-cancel}, \(F\) is a
  marked right fibration. By \myRef{prop-marked-rfib-pullback}, \(F'\)
  is a marked right fibration. By \myRef{prop-marked-rfib-cancel}, the
  composite \(S' \myXMorphism{F'} T' \myMorphism \myMaxMark C'\) is a
  marked right fibration. Hence \(S' \myMorphism C'\) is a cartesian
  fibration. By construction, \(F'\) is a cartesian functor.
\end{myProof}

\begin{myProposition}
  Any representable cartesian functor \(F \myElemOf S \myMorphism T\)
  over \(C\) is exponentiable in \(\myFiber{\myCart}{C}\). More
  precisely, \(F_{\myStar}\) is given by the pullback along
  \(F^{\myRepRA}\) which is indeed a functor
  \((\myFiber{\myCart}{C} \mySlice S) \myMorphism
  (\myFiber{\myCart}{C} \mySlice T)\) by
  \myRef{prop-cart-fib-pullback-by-ext}.
\end{myProposition}

\begin{myProof}
  By \myRef{prop-marked-rfib-cancel},
  \((\myFiber{\myCart}{C} \mySlice S)\) is equivalent to the category
  of marked right fibrations over \(S\). Then, the unit and counit for
  \(F^{\myStar} \myAdjRel F_{\myStar}\) are constructed by
  \myRef{cst-marked-rfib-2-functorial}.
\end{myProof}

\begin{myConstruction}
  Let \(C\) be a category with finite limits. Let
  \(\myRepClf_{C} \myDefEq \myCell_{1} \myPower C\) and we regard it
  as a cartesian fibration over \(C\) with the codomain
  projection. Its cartesian lift is given by
  pullback. \(\myRepClf_{C}\) has finite cartesian limits. Let
  \(\myRepClfPt_{C} \myDefEq (\myTerminal \mySlice_{C}
  \myRepClf_{C})\) and let
  \(\myGenRep_{C} \myElemOf \myRepClfPt_{C} \myMorphism
  \myRepClf_{C}\) denote the projection.
\end{myConstruction}

The cartesian functor \(\myGenRep_{C}\) is representable. Concretely,
\(\myRepClfPt_{C}\) is the category of diagrams in \(C\) of the form
\(
\begin{tikzcd}
  & y
  \arrow[d]
  \\ x
  \arrow[ur]
  \arrow[r, equal]
  & x
\end{tikzcd}
\),
and the right adjoint \(\myGenRepRA_{C}\) sends an object \(y
\myMorphism x\) in \(\myRepClf_{C}\) to the diagram
\(
\begin{tikzcd}
  & y \myBinProd_{x} y
  \arrow[d, "\myProj_{1}"]
  \\ y
  \arrow[ur, "\myDiagonal"]
  \arrow[r, equal]
  & y
\end{tikzcd}
\). One can see that \(\myGenRep_{C}\) is the universal representable
cartesian functor over \(C\), that is, every representable cartesian
functor over \(C\) is a pullback of \(\myGenRep_{C}\) in a unique way
(see \myCite[Section 5.1]{nguyen2022type-arxiv} for an analogous
result for representable maps of right fibrations). It then follows
that \(\myGenRep_{C}\) is a lex universe in
\(\myFiber{\myCart}{C}\). We have a canonical marked functor
\(\myWeakening_{T} \myElemOf T \myBinProd_{C} \myRepClf_{C}
\myMorphism \myApp{\myFam_{\myGenRep_{C}}}{T}\) over \(\myRepClf_{C}\)
corresponding to the projection
\((T \myBinProd_{C} \myRepClf_{C}) \myBinProd_{\myRepClf_{C}}
\myRepClfPt_{C} \myMorphism T\), where a morphism is marked if it is
cartesian over \(C\).

\begin{myDefinition}
  Let \(C\) be a category with finite limits and let \(T\) be a
  cartesian fibration over \(C\). We say \(T\) has
  \myDefine{\(C\)-coproducts} if
  \(\myWeakening_{T} \myElemOf T \myBinProd_{C} \myRepClf_{C}
  \myMorphism \myApp{\myFam_{\myGenRep_{C}}}{T}\) has a left adjoint
  \(\myCoprod\) over \(\myRepClf_{C}\) preserving marked morphisms.
\end{myDefinition}

The fiber of \(\myWeakening_{T}\) over an object \((f \myElemOf y
\myMorphism x) \myElemOf \myRepClf_{C}\) is the cartesian lift functor
\(f^{\myStar} \myElemOf \myFiber{T}{x} \myMorphism
\myFiber{T}{y}\). One can see that \(T\) has \(C\)-coproducts if and
only if every \(f^{\myStar} \myElemOf \myFiber{T}{x} \myMorphism
\myFiber{T}{y}\) has a left adjoint \(f_{\myBang}\) and the
Beck-Chevalley condition holds: for any pullback square
\[
  \begin{tikzcd}
    y'
    \arrow[r, "h"]
    \arrow[d, "f'"']
    & y
    \arrow[d, "f"]
    \\ x'
    \arrow[r, "g"']
    & x
  \end{tikzcd}
\]
in \(C\), the canonical natural transformation
\(f'_{\myBang} \myComp h^{\myStar} \myMMorphism g^{\myStar} \myComp
f_{\myBang} \myElemOf \myFiber{T}{y} \myMorphism \myFiber{T}{x'}\) is
an equivalence.

We have not found a standard terminology for the following concept.

\begin{myDefinition}
  Let
  \[
    \begin{tikzcd}
      T
      \arrow[r, "H"]
      \arrow[d]
      & T'
      \arrow[d]
      \\ C
      \arrow[r, "F"']
      & C'
    \end{tikzcd}
  \]
  be a morphism in \(\myCart\) from \(T \myMorphism C\) to
  \(T' \myMorphism C'\). Suppose that \(F\) has a left adjoint
  \(G\). We define a cartesian functor
  \(\myLeftComparison{H} \myElemOf \myApp{G^{\myStar}}{T} \myMorphism
  T'\) over \(C'\) by
  \(\myApp{\myLeftComparison{H}}{x \myComma t} \myDefEq
  \myApp{\eta^{\myStar}}{\myApp{H}{t}}\) where
  \(\eta \myElemOf x \myMorphism \myApp{F}{\myApp{G}{x}}\) is the unit of
  the adjunction \(G \myAdjRel F\). We say \((F \myComma H)\) is
  \myDefine{left-cartesian} if \(\myLeftComparison{H}\) is an
  equivalence.
\end{myDefinition}

By straightening-unstraightening, the functor \(H\) corresponds to a
natural transformation
\(T \myMMorphism T' \myComp \myApp{\myOp}{F} \myElemOf
\myApp{\myOp}{C} \myMorphism \myCat\). Then \((F \myComma H)\) is
left-cartesian if and only if the corresponding natural transformation
\(T \myComp \myApp{\myOp}{G} \myMMorphism T' \myElemOf
\myApp{\myOp}{C'} \myMorphism \myCat\) is an equivalence.

We end this subsection with a few facts about cartesian fibrations.

\begin{myProposition}
  \label{prop-pres-cart-fiber-pres}
  Let \(P \myElemOf T \myMorphism C\) be a morphism in \(\myPrR\). If
  \(P\) is a cartesian fibration, then every fiber \(\myFiber{T}{x}\)
  is presentable.
\end{myProposition}

\begin{myProof}
  The functor \(x \myElemOf \myTerminal \myMorphism C\) preserves
  filtered colimits and thus is accessible. Hence, \(\myFiber{T}{x}\)
  is the fiber product in \(\myAcc\) of \(P\) and \(x\). For any small
  diagram \(t \myElemOf I \myMorphism \myFiber{T}{x}\), the limit of
  \(t\) in \(\myFiber{T}{x}\) exists and is computed by the cartesian
  lift of the limit \(\myLim_{i \myElemOf I} t_{i}\) in \(T\) along
  the diagonal morphism \(x \myMorphism \myLim_{i \myElemOf I} x\). By
  \myRef{prop-acc-lim-pres}, \(\myFiber{T}{x}\) is presentable.
\end{myProof}

\begin{myProposition}
  \label{prop-cart-left-adj-unit-eqv}
  Let \(P \myElemOf T \myMorphism C\) be a cartesian fibration. Then
  \(P\) has a left adjoint with invertible unit if and only if every
  fiber \(\myFiber{T}{x}\) has an initial object. When this is the
  case, the left adjoint assigns the initial object in
  \(\myFiber{T}{x}\) to every \(x \myElemOf C\).
\end{myProposition}

\begin{myProof}
  Suppose that \(P\) has a left adjoint \(F\) with invertible
  unit. For every \(x \myElemOf C\) and
  \(t \myElemOf \myFiber{T}{x}\), the hom type
  \(\myApp{\myHom_{\myFiber{T}{x}}}{\myApp{F}{x} \myComma t}\) is the
  fiber of the equivalence
  \(P \myElemOf \myApp{\myHom_{T}}{\myApp{F}{x} \myComma t} \myEquiv
  \myApp{\myHom_{C}}{x \myComma x}\) over \(\myIdFun_{x}\) and thus
  contractible. Therefore, \(\myApp{F}{x}\) is an initial object in
  \(\myFiber{T}{x}\).

  Conversely, suppose that \(\myApp{F}{x}\) is an initial object in
  \(\myFiber{T}{x}\) for every \(x \myElemOf C\). Let
  \(t \myElemOf T\) and \(u \myElemOf x \myMorphism
  \myApp{P}{t}\). Take a cartesian morphism
  \(f \myElemOf t' \myMorphism t\) over \(u\). Then the fiber of
  \(P \myElemOf \myApp{\myHom_{T}}{\myApp{F}{x} \myComma t}
  \myMorphism \myApp{\myHom_{C}}{x \myComma \myApp{P}{t}}\) over \(u\)
  is equivalent to
  \(\myApp{\myHom_{\myFiber{T}{x}}}{\myApp{F}{x} \myComma t'}\) which
  is contractible as \(\myApp{F}{x}\) is initial. Therefore, \(F\)
  defines a left adjoint of \(P\) with invertible unit.
\end{myProof}

\begin{myCorollary}
  \label{prop-pres-cart-left-adj-ff}
  Let \(P \myElemOf T \myMorphism C\) be a morphism in \(\myPrR\). If
  \(P\) is a cartesian fibration, then the unit of the left adjoint
  \(Q\) of \(P\) is an equivalence, and \(\myApp{Q}{x}\) is the
  initial object in the fiber \(\myFiber{T}{x}\) for every
  \(x \myElemOf C\). \myQED
\end{myCorollary}

\section{Higher inductive types in \(\myType\)}
\label{sec-limit-theories}

We consider higher inductive types in the category \(\myType\) of
small types. Higher inductive types are understood as \emph{initial
  algebras} for some kind of algebraic theories
\myCite{awodey2012inductive,sojakova2015higher,awodey2017homotopy}. The
notion of algebraic theory suitable in the context of dependent type
theory is generalized algebraic theory
\myCite{cartmell1978generalised}, which has the same expressive power
as essentially algebraic theory \myCite{freyd1972aspects}, cartesian
theory \myCite{johnstone2002sketches2}, partial Horn theory
\myCite{palmgren2007partial}, and limit theory
\myCite{coste1979localisation}. All of them are identified with
categories with finite limits via the construction of syntactic
categories.

\begin{myConstruction}
  We define the category \(\myFinLimTh\) of small \myDefine{(finite)
    limit theories (in \(\myType\))} to be
  \(\myApp{\myCoreI}{\myLex}\).
\end{myConstruction}

For a limit theory \(t\), the category of small \myDefine{algebras for
  \(t\)} is the full subcategory
\(t \myLexPower \myType \mySub t \myPower \myType\) consisting of
those functors preserving finite limits. The functor
\(t \myMapsTo t \myLexPower \myType \myElemOf
\myApp{\myOp}{\myFinLimTh} \myMorphism \myEnlarge \myCat\) is in fact
induced by a universal property of \(\myApp{\myOp}{\myFinLimTh}\)
(\myRef{prop-lim-th-initial} and
\myRef{prop-algebra-in-type-lex-functor}). Consider the following
structures on \(\myApp{\myOp}{\myFinLimTh}\); see \myCite[Section
5.3]{nguyen2022type-arxiv} for more details.
\begin{itemize}
\item \(\myApp{\myOp}{\myFinLimTh}\) has small limits.
\item \(\myApp{\myOp}{\myFinLimTh}\) is equipped with a lex universe
  \(\myUniv \myElemOf \myUnivTot \myMorphism
  \myUnivBase\). Concretely, \(\myUnivBase\) in
  \(\myApp{\myCoreI}{\myLex}\) is the category with finite limits
  freely generated by an object \(w\), and \(\myUnivTot\) is the
  extension of \(\myUnivBase\) by a morphism
  \(\myTerminal \myMorphism w\).
\end{itemize}

\begin{myProposition}[{\myCite[Corollary 5.21]{nguyen2022type-arxiv}}]
  \label{prop-lim-th-initial}
  \(\myApp{\myOp}{\myFinLimTh}\) is the initial category with small
  limits and a lex universe. \myQED
\end{myProposition}

We give an explicit description of the unique morphism
\(\myApp{\myOp}{\myFinLimTh} \myMorphism C\) of categories with small
limits and a lex universe (\myRef{prop-unique-morphism-lex-univ}).

\begin{myLemma}
  \label{prop-lex-slice}
  Let \(C \myElemOf \myFinLimTh\) and
  \(x \myElemOf \myUnivBase \myMorphism C\). We regard \(x\) as an
  object in \(C\). Then the slice \((C \mySlice x)\) fits into the
  following pushout square in \(\myFinLimTh\).
  \[
    \begin{tikzcd}
      \myUnivBase
      \arrow[r, "x"]
      \arrow[d, "\myUniv"']
      & C
      \arrow[d]
      \\ \myUnivTot
      \arrow[r]
      & (C \mySlice x)
    \end{tikzcd}
  \]
\end{myLemma}

\begin{myProof}
  This follows from \myCite[Proposition 3.25]{nguyen2022type-arxiv}.
\end{myProof}

\begin{myConstruction}
  Let \(C\) be a category with small limits and a lex universe
  \(\myUniv_{C}\). For an object \(x \myElemOf C\), we define
  \((C \mySlice^{\myUniv_{C}} x) \mySub (C \mySlice x)\) to be the
  full subcategory consisting of \(\myUniv_{C}\)-small morphisms to
  \(x\). Note that
  \((C \mySlice^{\myUniv_{C}} x) \mySub (C \mySlice x)\) is closed
  under finite limits.
\end{myConstruction}

\begin{myProposition}
  \label{prop-unique-morphism-lex-univ}
  Let \(C\) be a category with small limits and a lex universe
  \(\myUniv_{C} \myElemOf \myUnivTot_{C} \myMorphism
  \myUnivBase_{C}\). Suppose that \(C\) is in
  \(\myEnlarge^{n} \myCat\) for \(n \myGe 1\). Let
  \(F \myElemOf \myApp{\myOp}{\myFinLimTh} \myMorphism C\) denote the
  unique morphism of categories with small limits and a lex
  universe. Then we have an equivalence
  \begin{equation*}
    \myApp{\myHom_{C}}{x \myComma \myApp{F}{t}} \myEquiv
    \myApp{\myHom_{\myEnlarge^{n} \myFinLimTh}}{t \myComma (C
      \mySlice^{\myUniv_{C}} x)}
  \end{equation*}
  natural in \(x \myElemOf C\) and \(t \myElemOf \myFinLimTh\).
\end{myProposition}

\begin{myProof}
  The statement is equivalent to that the functor
  \[
    H \myDefEq (t \myMapsTo (x \myMapsTo \myApp{\myHom}{t \myComma (C
      \mySlice^{\myUniv_{C}} x)})) \myElemOf
    \myApp{\myOp}{\myFinLimTh} \myMorphism \myApp{\myPsh}{C}
  \]
  is equivalent to
  \(F \myElemOf \myApp{\myOp}{\myFinLimTh} \myMorphism C \mySub
  \myApp{\myPsh}{C}\). Because the Yoneda embedding preserves all
  existing limits and exponentials,
  \(F \myElemOf \myApp{\myOp}{\myFinLimTh} \myMorphism
  \myApp{\myPsh}{C}\) is a morphism of categories with small limits
  and a lex universe, where \(\myUniv_{C}\) in
  \(C \mySub \myApp{\myPsh}{C}\) is chosen as the lex universe in
  \(\myApp{\myPsh}{C}\). By the initiality of
  \(\myApp{\myOp}{\myFinLimTh}\), it suffices to show that \(H\) is a
  morphism of categories with small limits and a lex universe. Clearly
  \(H\) preserves small limits. \(H\) takes \(\myUnivBase\) to
  \(\myUnivBase_{C}\) as
  \(\myApp{\myHom}{\myUnivBase \myComma (C \mySlice^{\myUniv_{C}} x)}
  \myEquiv \myApp{\myObj}{C \mySlice^{\myUniv_{C}} x} \myEquiv
  \myApp{\myHom}{x \myComma \myUnivBase_{C}}\) and similarly
  \(\myUniv\) to \(\myUniv_{C}\). For any
  \(A \myElemOf \myApp{\myPsh}{C}\), we have
  \(\myApp{\myApp{\myFam_{\myUniv_{C}}}{A}}{x} \myEquiv (y \myElemOf
  \myApp{\myHom}{x \myComma \myUnivBase_{C}}) \myBinProd
  \myApp{A}{\myApp{y^{\myStar}}{\myUnivTot_{C}}}\). Then
  \begin{align*}
    & \myApp{\myHom}{\myApp{\myFam_{\myUniv}}{t} \myComma (C
    \mySlice^{\myUniv_{C}} x)}
    \\ \myEquiv{}
    & (y \myElemOf \myApp{\myHom}{\myUnivBase \myComma (C
      \mySlice^{\myUniv_{C}} x)}) \myBinProd \myApp{\myHom}{t \myComma
      (C \mySlice^{\myUniv_{C}} x) \myBinCoprod_{\myUnivBase}
      \myUnivTot}
    \\ \myEquiv{}
    & (y \myElemOf \myApp{\myObj}{C \mySlice^{\myUniv_{C}} x})
      \myBinProd \myApp{\myHom}{t \myComma ((C \mySlice^{\myUniv_{C}}
      x) \mySlice y)}
      \tag{\myRef{prop-lex-slice}}
    \\ \myEquiv{}
    & (y \myElemOf \myApp{\myObj}{C \mySlice^{\myUniv_{C}} x})
      \myBinProd \myApp{\myHom}{t \myComma (C \mySlice^{\myUniv_{C}}
      y)}
      \tag{\myRef{prop-small-morphisms-cancel}}
    \\ \myEquiv{}
    & (y \myElemOf \myApp{\myHom}{x \myComma \myUnivBase_{C}})
      \myBinProd \myApp{\myHom}{t \myComma (C \mySlice^{\myUniv_{C}}
      \myApp{y^{\myStar}}{\myUnivTot_{C}})}
  \end{align*}
  and thus
  \(\myApp{H}{\myApp{\myFam_{\myUniv}}{t}} \myEquiv
  \myApp{\myFam_{\myUniv_{C}}}{\myApp{H}{t}}\).
\end{myProof}

Recall that the projection
\(\myUnivLFib \myElemOf (\myTerminal \mySlice \myType) \myMorphism
\myType\) is the universal left fibration with small fibers.

\begin{myProposition}
  \label{prop-type-lex-universe-in-cat}
  The projection
  \(\myUnivLFib \myElemOf (\myTerminal \mySlice \myType) \myMorphism
  \myType\) is a lex universe in \(\myEnlarge \myCat\).
\end{myProposition}

\begin{myProof}
  \(\myUnivLFib\) is exponentiable by \myCite[Lemma
  3.2.1]{ayala2020fibrations}. The other conditions follow from the
  universal property.
\end{myProof}

By \myRef{prop-lim-th-initial} and
\myRef{prop-type-lex-universe-in-cat}, we have a unique morphism
\(\myFinLimThAlg \myElemOf \myApp{\myOp}{\myFinLimTh} \myMorphism
\myEnlarge \myCat\) of categories with small limits and a lex
universe, which we call the \myDefine{semantics
  functor}.

\begin{myProposition}
  \label{prop-algebra-in-type-lex-functor}
  \(\myApp{\myFinLimThAlg}{t} \myEquiv t \myLexPower \myType\)
  naturally in \(t \myElemOf \myFinLimTh\).
\end{myProposition}

\begin{myProof}
  We have
  \(\myApp{\myHom_{\myEnlarge \myCat}}{C \myComma t \myLexPower
    \myType} \myEquiv \myApp{\myHom_{\myEnlarge \myFinLimTh}}{t
    \myComma C \myPower \myType} \myEquiv \myApp{\myHom_{\myEnlarge^{2}
      \myFinLimTh}}{t \myComma (\myEnlarge \myCat
    \mySlice^{\myUnivLFib} C)}\). Then
  \(\myApp{\myFinLimThAlg}{t} \myEquiv t \myLexPower \myType\) by
  \myRef{prop-unique-morphism-lex-univ} and by Yoneda.
\end{myProof}

\begin{myCorollary}
  \label{prop-initial-algebras-in-type}
  \(\myApp{\myFinLimThAlg}{t}\) has an initial object for every
  \(t \myElemOf \myFinLimTh\). \myQED
\end{myCorollary}

\section{Categories with higher inductive types}
\label{sec-inductive-types}

We consider higher inductive types in a general category \(C\) with
finite limits. The idea is to perform the construction of
\(\myFinLimTh\) and \(\myFinLimThAlg\) explained in
\myRef{sec-limit-theories} \emph{over the base category \(C\)} instead
of \(\myType\). Initial algebras for limit theories over \(C\) may not
exist. We say \(C\) has higher inductive types if every limit theory
over \(C\) has an initial algebra. We begin by introducing a category
where the category of limit theories over \(C\) live.

\begin{myDefinition}
  A \myDefine{(finite) limit metatheory} consists of the following
  data.
  \begin{itemize}
  \item A category \(C\) with finite limits.
  \item A cartesian fibration \(T\) over \(C\) that has finite
    cartesian colimits and \(C\)-coproducts.
  \item A cartesian co-lex-universe
    \(\myUniv_{T} \myElemOf \myUnivBase_{T} \myMorphism
    \myUnivTot_{T}\) in \(T\).
  \end{itemize}
  We have the \(\myFinitary\)-presentable \(2\)-category
  \(\myFinLimMeta\) of small limit metatheories.
\end{myDefinition}

\begin{myRemark}
  \label{rem-limit-metatheory-and-signatures}
  The notion of limit metatheory is almost a categorical counterpart
  of the theory of signatures introduced in
  \myCite{kaposi2020signatures}. We summarize the correspondence in
  \myRef{tab-limit-metatheory-and-signature}. The universe in a limit
  metatheory is required to be univalent, while the universe in the
  theory of signatures is not. This is not an essential difference
  because univalence is not needed for defining signatures for higher
  inductive types. We impose univalence for better categorical
  behavior. A structure corresponding to infinitary parameters is not
  considered in the present paper.
  \begin{table}
    \begin{tabular}{ll}
      Structure of \(\myApp{\myOp_{C}}{T}\)
      & Rule group for signatures
      \\ \hline
      \(\myUniv_{T}\)
      & (2) Universe
      \\ Univalence of \(\myUniv_{T}\)
      & (None)
      \\ Exponentials by \(\myUniv_{T}\)
      & (3) Inductive parameters
      \\ Identity types in \(\myUniv_{T}\)
      & (4) Paths between point and path constructors
      \\ Finite cartesian limits
      & (5) Paths between type constructors
      \\ \(C\)-products
      & (6) External parameters
      \\ (None)
      & (7) Infinitary parameters
    \end{tabular}
    \caption{Limit metatheory and the theory of signatures}
    \label{tab-limit-metatheory-and-signature}
    Each structure of \(\myApp{\myOp_{C}}{T}\) on the left corresponds
    to the rule group of the theory of signatures found in
    \myCite[Figures 1 and 2]{kaposi2020signatures} on the right in the
    same row.
  \end{table}
\end{myRemark}

\begin{myNotation}
  We write \(\myFinLimTh \myElemOf \myLex \myMorphism \myFinLimMeta\)
  for the left adjoint of the forgetful functor
  \(\myFinLimMeta \myMorphism \myLex\).
\end{myNotation}

\begin{myLemma}
  \label{prop-lim-meta-lex-cart}
  The forgetful functor \(\myFinLimMeta \myMorphism \myLex\) is a
  cartesian fibration.
\end{myLemma}

\begin{myProof}
  The cartesian lift is given by pullback.
\end{myProof}

\begin{myLemma}
  \label{prop-limth-initial-in-fiber}
  The unit of the left adjoint \(\myFinLimTh\) is invertible, and
  \(\myApp{\myFinLimTh}{C}\) is the initial object in the fiber
  \(\myFiber{\myFinLimMeta}{C}\) for every \(C \myElemOf \myLex\).
\end{myLemma}

\begin{myProof}
  By \myRef{prop-lim-meta-lex-cart} and
  \myRef{prop-pres-cart-left-adj-ff}.
\end{myProof}

The semantics functor introduced in \myRef{sec-limit-theories} is a
functor \(\myApp{\myOp}{\myFinLimTh} \myMorphism \myEnlarge
\myCat\). By straightening-unstraightening, it corresponds to a
cartesian fibration over \(\myFinLimTh\). We define a semantics of a
limit metatheory in the fibrational form.

\begin{myDefinition}
  \label{def-semantics}
  Let \((C \myComma T)\) be a limit metatheory. A \myDefine{semantics
    of \(T\)} is a cartesian fibration \(A\) over \(T\) satisfying the
  following conditions.
  \begin{enumerate}
  \item \label[condition]{cond-sem-fin-colim} For any finite simplicial type
    \(I\), the square
    \[
      \begin{tikzcd}
        A
        \arrow[r, "\myDiagonal"]
        \arrow[d]
        & I \myMkPower_{C} A
        \arrow[d]
        \\ T
        \arrow[r, "\myDiagonal"']
        & I \myPower_{C} T,
      \end{tikzcd}
    \]
    where the bottom functor has the left adjoint \(\myColim\), is
    left-cartesian. Here \(I \myMkPower_{C} A\) is the full
    subcategory of \(I \myPower_{C} A\) consisting of those diagram
    \(a \myElemOf I \myMorphism A\) such that
    \(a_{i} \myMorphism a_{j}\) is cartesian over \(T\) for every edge
    \(i \myMorphism j\) in \(I\).
  \item \label[condition]{cond-sem-coprod} The square
    \[
      \begin{tikzcd}
        A \myBinProd_{C} \myRepClf_{C}
        \arrow[r, "\myWeakening_{A}"]
        \arrow[d]
        & \myApp{\myFam_{\myGenRep_{C}}}{A}
        \arrow[d]
        \\ T \myBinProd_{C} \myRepClf_{C}
        \arrow[r, "\myWeakening_{T}"']
        & \myApp{\myFam_{\myGenRep_{C}}}{T},
      \end{tikzcd}
    \]
    where the bottom functor has the left adjoint \(\myCoprod\), is
    left-cartesian.
  \item \label[condition]{cond-sem-univ} The functor
    \(\myApp{\myUniv_{T}^{\myStar}}{A} \myElemOf
    \myApp{\myUnivTot_{T}^{\myStar}}{A} \myMorphism
    \myApp{\myUnivBase_{T}^{\myStar}}{A}\) is equivalent to
    \(\myGenRep_{C} \myElemOf \myRepClfPt_{C} \myMorphism
    \myRepClf_{C}\). Such an equivalence is unique because of the
    univalence of \(\myGenRep_{C}\).
  \item \label[condition]{cond-sem-pushforward} The square
    \[
      \begin{tikzcd}
        A \myBinProd_{T} (\myUnivBase_{T} \mySlice_{C} T)
        \arrow[r]
        \arrow[d]
        & \myApp{(\myGenRep_{C} \myBinProd_{C} (\myUnivTot_{T} \mySlice_{C} T))_{\myStar}}{A \myBinProd_{T} (\myUnivTot_{T} \mySlice_{C} T)}
        \arrow[d]
        \\ (\myUnivBase_{T} \mySlice_{C} T)
        \arrow[r, "(\myUniv_{T})_{\myBang}"']
        & (\myUnivTot_{T} \mySlice_{C} T),
      \end{tikzcd}
    \]
    where the bottom pushout functor has a left adjoint, is
    left-cartesian. Here, the top functor is constructed as
    follows. It suffices to construct a functor
    \(A \myBinProd_{T} (\myUnivBase_{T} \mySlice_{C} T) \myMorphism
    \myApp{(\myUniv_{T})_{\myBang}^{\myStar}}{\myApp{(\myGenRep_{C}
        \myBinProd_{C} (\myUnivTot_{T} \mySlice_{C} T))_{\myStar}}{A
        \myBinProd_{T} (\myUnivTot_{T} \mySlice_{C} T)}}\) over
    \((\myUnivBase_{T} \mySlice_{C} T)\). The codomain is equivalent
    to
    \(\myApp{(\myGenRep_{C} \myBinProd_{C} (\myUnivBase_{T}
      \mySlice_{C}
      T))_{\myStar}}{\myApp{(\myUniv_{T})_{\myBang}^{\myStar}}{A
        \myBinProd_{T} (\myUnivTot_{T} \mySlice_{C} T)}}\). This is
    further equivalent to
    \(\myApp{(\myGenRep_{C} \myBinProd_{C} (\myUnivBase_{T}
      \mySlice_{C} T))_{\myStar}}{\myApp{(\myGenRep_{C} \myBinProd_{C}
        (\myUnivBase_{T} \mySlice_{C} T))^{\myStar}}{A \myBinProd_{T}
        (\myUnivBase_{T} \mySlice_{C} T)}}\), because
    \myRef{cond-sem-fin-colim} and \myRef{cond-sem-univ}
    imply that the fiber of \(A\) over a pushout along \(\myUniv_{T}\)
    is equivalent to a pullback along \(\myGenRep_{C}\). Thus, the
    unit of the adjunction
    \((\myGenRep_{C} \myBinProd_{C} (\myUnivBase_{T} \mySlice_{C}
    T))^{\myStar} \myAdjRel (\myGenRep_{C} \myBinProd_{C}
    (\myUnivBase_{T} \mySlice_{C} T))_{\myStar}\) gives the required
    functor.
  \end{enumerate}
  Recall that \(P_{\myStar}\) and \(\myFam_{P}\) for a representable
  cartesian functor \(P\) over \(C\) are defined in terms of the
  pullback along \(P^{\myRepRA}\). We thus have the
  \(\myFinitary\)-presentable \(2\)-category \(\myFinLimSem\) of tuples
  \((C \myComma T \myComma A)\) where
  \((C \myComma T) \myElemOf \myFinLimMeta\) and \(A\) is a small
  semantics of \(T\).
\end{myDefinition}

With \myRef{prop-cart-left-adj-unit-eqv} in mind, initial algebras are
defined in terms of left adjoints with invertible units.

\begin{myDefinition}
  Let \((C \myComma T)\) be a limit metatheory and let \(A\) be a
  semantics of \(T\). We say \(A\) \myDefine{has initial algebras}
  if \(A \myMorphism T\) has a cartesian left adjoint over \(C\) with
  invertible unit. We have the \(\myFinitary\)-presentable \(2\)-category
  \(\myFinLimInit\) of tuples \((C \myComma T \myComma A)\) where
  \((C \myComma T \myComma A) \myElemOf \myFinLimSem\) and \(A\)
  has initial algebras.
\end{myDefinition}

In the rest of this section, we forget about the \(2\)-category
structures on \(\myLex\), \(\myFinLimMeta\), \(\myFinLimSem\), and
\(\myFinLimInit\) and regard them as \(\myFinitary\)-presentable
categories to avoid diving too much into \(2\)-category theory. The
category \(\myFinIndType\) constructed below may be extended to an
\(\myFinitary\)-presentable \(2\)-category, but the \(1\)-category
structure is enough to state and prove canonicity
(\myRef{sec-canonicity}). Still, we will use the \(2\)-category
structure on \(\myFinLimSem\) in \myRef{sec-canonicity}.

\begin{myConstruction}
  We define the category \(\myFinIndType\) of small categories with
  \myDefine{(finitary) higher inductive types} to be the fiber product
  of \(\myFinLimTh \myElemOf \myLex \myMorphism \myFinLimMeta\) and
  the forgetful functor \(\myFinLimInit \myMorphism \myFinLimMeta\).
\end{myConstruction}

\(\myFinIndType\) fits into the following pullback diagram in
\(\myEnlarge \myCat\).
\begin{equation}
  \label[diagram]{eq-ind-type-pb}
  \begin{tikzcd}
    \myFinIndType
    \arrow[r]
    \arrow[d]
    & \myFinLimInit
    \arrow[d]
    \\ \myApp{\myFinLimTh^{\myStar}}{\myFinLimSem}
    \arrow[r]
    \arrow[d]
    & \myFinLimSem
    \arrow[d]
    \\ \myLex
    \arrow[r, "\myFinLimTh"']
    & \myFinLimMeta
  \end{tikzcd}
\end{equation}
Unwinding the definitions, an object in \(\myFinIndType\) is a pair
\((C \myComma A)\) where \(C \myElemOf \myLex\) and \(A\) is a
semantics of \(\myApp{\myFinLimTh}{C}\) that has initial
algebras. We will see in \myRef{prop-lim-sem-lex-eqv} below that there
exists a unique semantics of \(\myApp{\myFinLimTh}{C}\) which we refer
to as \(\myApp{\myFinLimThAlg}{C}\). Thus, the component \(A\) is
redundant, and a category \(C\) with finite limits is in
\(\myFinIndType\) if and only if \(\myApp{\myFinLimThAlg}{C}\) has
initial algebras.

In the rest of this section, we prove that
\(\myFinIndType \myMorphism \myLex\) is in \(\myPrRFin\) and monic
(\myRef{prop-ind-type-pres}). Note that, because
\(\myFinLimTh \myElemOf \myLex \myMorphism \myFinLimMeta\) is
\emph{not} in \(\myPrRFin\), \myRef{eq-ind-type-pb} is
not pullbacks in \(\myPrRFin\), and thus it is not immediate that
\(\myFinIndType\) is \(\myFinitary\)-presentable.

\begin{myConstruction}
  Let \(C\) be a category with finite limits. We define a cartesian
  fibration \(\myCat_{C}\) over \(C\) to be the pullback of
  \(\myCart \myMorphism \myCat\) along the functor
  \(x \myMapsTo (C \mySlice x) \myElemOf C \myMorphism \myCat\) whose
  action on morphisms is postcomposition. \(\myCat_{C}\) has
  \(C\)-products. Indeed, the cartesian lift functor
  \(f^{\myStar} \myElemOf \myFiber{\myCat_{C}}{x} \myMorphism
  \myFiber{\myCat_{C}}{y}\) for \(f \myElemOf y \myMorphism x\) is
  equivalent to the precomposition functor
  \(\myApp{\myOp}{C \mySlice x} \myPower \myCat \myMorphism
  \myApp{\myOp}{C \mySlice y} \myPower \myCat\) which has a right
  adjoint. The Beck-Chevalley condition follows from the base change
  formula \myCite[Theorem 6.4.13]{cisinski2019higher} since
  \((C \mySlice x') \myMorphism (C \mySlice x)\) for
  \(g \myElemOf x' \myMorphism x\) is a right fibration and thus
  smooth. Each fiber
  \(\myFiber{\myCat_{C}}{x} \myEquiv \myFiber{\myCart}{C \mySlice x}\)
  has the lex universe
  \(\myGenRep_{(C \mySlice x)} \myElemOf \myRepClfPt_{(C \mySlice x)}
  \myMorphism \myRepClf_{(C \mySlice x)}\). This defines a cartesian
  lex universe in \(\myCat_{C}\) by which we regard
  \(\myApp{\myOp_{C}}{\myCat_{C}}\) as a limit metatheory over \(C\).
\end{myConstruction}

\(\myCat_{C}\) classifies cartesian fibrations as
\begin{align*}
  & \myApp{\myHom_{\myFiber{\myEnlarge \myCart}{C}}}{\myApp{\myOp_{C}}{T} \myComma
    \myCat_{C}}
  \\ \myEquiv {}
  & \myEnd_{x \myElemOf C} \myApp{\myHom_{\myEnlarge
    \myCat}}{\myApp{\myOp}{\myFiber{T}{x}} \myComma
    \myFiber{\myCart}{(C \mySlice x)}}
  \\ \myEquiv {}
  & \myEnd_{x \myElemOf C} \myApp{\myHom_{\myEnlarge
    \myCat}}{\myApp{\myOp}{\myFiber{T}{x}} \myComma \myApp{\myOp}{C
    \mySlice x} \myPower \myCat}
  \\ \myEquiv {}
  & \myApp{\myHom_{\myEnlarge \myCat}}{\myApp{\myOp}{\myEnd^{x \myElemOf C} (C
    \mySlice x) \myBinProd \myFiber{T}{x}} \myComma
    \myCat}
  \\ \myEquiv {}
  & \myApp{\myHom_{\myEnlarge \myCat}}{\myApp{\myOp}{T} \myComma
    \myCat}
  \\ \myEquiv {}
  & \myFiber{\myCart}{T}
\end{align*}
for every cartesian fibration \(T\) over \(C\). This is a special case
of internal straightening-unstraightening
\myCite{martini2022cocartesian-arxiv}.

\begin{myLemma}
  \label{prop-semantics-straightening}
  Let \((C \myComma T) \myElemOf \myFinLimMeta\). Then a cartesian
  fibration \(A\) over \(T\) is a semantics of \(T\) if and only if
  the corresponding cartesian functor
  \(T \myMorphism \myApp{\myOp_{C}}{\myCat_{C}}\) over \(C\) is a
  morphism of limit metatheories.
\end{myLemma}

\begin{myProof}
  \myRef{def-semantics} is designed to satisfy this property.
\end{myProof}

\begin{myLemma}
  \label{prop-lim-sem-lex-eqv}
  \(\myApp{\myFinLimTh^{\myStar}}{\myFinLimSem} \myEquiv \myLex\).
\end{myLemma}

\begin{myProof}
  We first show that the object part of the forgetful functor
  \(U \myElemOf \myApp{\myFinLimTh^{\myStar}}{\myFinLimSem}
  \myMorphism \myLex\) is an equivalence. Let \(C \myElemOf
  \myLex\). By \myRef{prop-semantics-straightening}, the type of
  semantics of \(\myApp{\myFinLimTh}{C}\) is equivalent to the type of
  morphisms
  \(\myApp{\myFinLimTh}{C} \myMorphism \myApp{\myOp_{C}}{\myCat_{C}}\)
  of limit metatheories over \(C\). The latter is contractible by
  \myRef{prop-limth-initial-in-fiber}. Therefore, the object part of
  \(U\) is an equivalence.

  For the morphism part of \(U\), repeat the whole argument internally
  to \(\myCell_{1} \myPower \myEnlarge \myType\) instead of
  \(\myEnlarge \myType\). Let
  \(U_{\myCell_{1}} \myElemOf
  \myApp{\myFinLimTh_{\myCell_{1}}^{\myStar}}{\myFinLimSem_{\myCell_{1}}}
  \myMorphism \myLex_{\myCell_{1}}\) be the functor constructed in the
  same way as \(U\) but internally to
  \(\myCell_{1} \myPower \myEnlarge \myType\). Then
  \(\myCell_{1} \myPower U\) is the image of \(U_{\myCell_{1}}\) by
  the domain projection
  \(\myCell_{1} \myPower \myEnlarge \myType \myMorphism \myEnlarge
  \myType\). The object part of \(U_{\myCell_{1}}\) is an equivalence
  (internally to \(\myCell_{1} \myPower \myEnlarge \myType\)) as we
  have seen. The object part of \(\myCell_{1} \myPower U\) is thus an
  equivalence (internally to \(\myEnlarge \myType\)), but it is the
  morphism part of \(U\).
\end{myProof}

\begin{myNotation}
  We refer to the composite
  \(\myLex \myEquiv \myApp{\myFinLimTh^{\myStar}}{\myFinLimSem}
  \myMorphism \myFinLimSem\) as \(\myFinLimThAlg\).
\end{myNotation}

\begin{myConvention}
  \label{conv-marked-lim-meta}
  We regard \(\myFinLimMeta\), \(\myFinLimSem\), and
  \(\myFinLimInit\) as marked categories as follows.
  \begin{itemize}
  \item A morphism
    \((F \myComma G) \myElemOf (C \myComma T) \myMorphism (C' \myComma
    T')\) in \(\myFinLimMeta\) is marked if
    \(F \myElemOf C \myMorphism C'\) is an equivalence.
  \item A morphism
    \((F \myComma G \myComma H) \myElemOf (C \myComma T \myComma A)
    \myMorphism (C' \myComma T' \myComma A')\) in \(\myFinLimSem\) or
    \(\myFinLimInit\) is marked if \((F \myComma G)\) in
    \(\myFinLimMeta\) is marked and
    \((G \myComma H) \myElemOf (T \myComma A) \myMorphism (T' \myComma
    A')\) is a pullback.
  \end{itemize}
\end{myConvention}

\begin{myLemma}
  \label{prop-lim-meta-marked-rfib}
  Under \myRef{conv-marked-lim-meta}, the forgetful functors
  \(\myFinLimInit \myMorphism \myFinLimSem\) and
  \(\myFinLimSem \myMorphism \myFinLimMeta\) are marked right
  fibrations. \myQED
\end{myLemma}

\begin{myLemma}
  \label{prop-ind-type-corefl-in-lim-meta-init}
  The functor \(\myFinIndType \myMorphism \myFinLimInit\) is fully
  faithful and has a right adjoint, and the mate of
  \myRef{eq-ind-type-pb} is an equivalence.
\end{myLemma}

\begin{myProof}
  By \myRef{prop-limth-initial-in-fiber} \(\myFinLimTh\) is fully
  faithful, and so is its pullback
  \(\myFinIndType \myMorphism \myFinLimInit\). It also follows that
  the counit of \(\myFinLimTh\) is component-wise marked under
  \myRef{conv-marked-lim-meta}. Thus, by
  \myRef{prop-lift-marked-adjoint} and
  \myRef{prop-lim-meta-marked-rfib},
  \(\myFinIndType \myMorphism \myFinLimInit\) has a right adjoint and
  the mate of \myRef{eq-ind-type-pb} is an equivalence.
\end{myProof}

\begin{myTheorem}
  \label{prop-ind-type-pres}
  \(\myFinIndType\) is \(\myFinitary\)-presentable, and the forgetful
  functor \(\myFinIndType \myMorphism \myLex\) is in \(\myPrRFin\) and
  monic.
\end{myTheorem}

\begin{myProof}
  \myRef{eq-ind-type-pb} is a pullback diagram in
  \(\myAccFin\). By \myRef{prop-ind-type-corefl-in-lim-meta-init},
  \(\myFinIndType\) has small limits computed by the coreflection of
  small limits in \(\myFinLimInit\). Then the forgetful functor
  \(\myFinIndType \myMorphism \myLex\) preserves small limits. By
  \myRef{prop-acc-lim-pres}, \(\myFinIndType\) is presentable. By
  \myRef{prop-lim-sem-lex-eqv} and because
  \(\myFinLimInit \myMorphism \myFinLimSem\) is monic,
  \(\myFinIndType \myMorphism \myLex\) is monic. It is also
  conservative by \myRef{prop-mono-implies-conservative}. Then, by
  \myRef{prop-fin-pres-closed-under-csv}, \(\myFinIndType\) is
  \(\myFinitary\)-presentable.
\end{myProof}

\begin{myCorollary}
  \label{prop-when-lex-has-ind-type}
  A category \(C\) with finite limits has higher inductive types if
  and only if there exist a limit metatheory \(T\) over \(C\) and a
  semantics \(A\) of \(T\) that has initial algebras.
\end{myCorollary}

\begin{myProof}
  The ``only if'' part is trivial. For the ``if'' part, apply the
  coreflection obtained in
  \myRef{prop-ind-type-corefl-in-lim-meta-init}.
\end{myProof}

\section{Examples of higher inductive types}
\label{sec-examples-of-higher-inductive-types}

We give examples of higher inductive types: finite colimits
(\myRef{prop-ind-type-fin-colim}); natural number object
(\myRef{prop-ind-type-natural-numbers}); truncation
(\myRef{prop-ind-type-truncation}).

\begin{myRemark}
  Having the relationship between limit metatheories and the theory of
  signatures explained in \myRef{rem-limit-metatheory-and-signatures}
  in mind, one can construct more complex examples of higher inductive
  types such as indexed inductive types, mutual inductive types, and
  higher inductive-inductive types. Presenting such higher inductive
  types in categorical terms is, however, extremely complicated. An
  interpretation of the theory of signatures in a limit metatheory
  should be established, but we leave it as future work.
\end{myRemark}

\begin{myDefinition}
  We say a structure on a category \(C\) with finite limits is
  \myDefine{stable under pullback} if every slice \((C \mySlice x)\)
  has that structure and if the pullback functor
  \(f^{\myStar} \myElemOf (C \mySlice x) \myMorphism (C \mySlice x')\)
  preserves that structure for every morphism
  \(f \myElemOf x' \myMorphism x\) in \(C\).
\end{myDefinition}

In the following, we claim that some structure is stable under
pullback and preserved by any morphism in \(\myFinIndType\) but omit
its proof, because it is clear from the construction of that structure
from higher inductive types.

\begin{myConstruction}
  Let \((C \myComma T) \myElemOf \myFinLimMeta\). We define
  \((a \myComma t) \myMapsTo a \myAct t \myElemOf \myRepClf_{C}
  \myBinProd_{C} T \myMorphism T\) to be the composite
  \(\myRepClf_{C} \myBinProd_{C} T \myXMorphism{\myWeakening}
  \myApp{\myFam_{\myGenRep_{C}}}{T} \myXMorphism{\myCoprod}
  \myRepClf_{C} \myBinProd_{C} T \myXMorphism{\myProj_{2}} T\). We
  also define
  \((a \myMapsTo a \myAct_{\myUnivBase_{T}} \myUnivTot_{T}) \myElemOf
  \myRepClf_{C} \myMorphism (\myUnivBase_{T} \mySlice_{C} T)\) by the
  pushout
  \[
    \begin{tikzcd}
      a \myAct \myUnivBase_{T}
      \arrow[r, "\epsilon"]
      \arrow[d]
      & \myUnivBase_{T}
      \arrow[d]
      \\ a \myAct \myUnivTot_{T}
      \arrow[r]
      & a \myAct_{\myUnivBase_{T}} \myUnivTot_{T}
    \end{tikzcd}
  \]
  where \(\epsilon\) is derived from the counit of the adjunction
  \(\myCoprod \myAdjRel \myWeakening\).
\end{myConstruction}

\begin{myProposition}
  \label{prop-ind-type-fin-colim}
  Let \(C\) be a category with higher inductive types. Then \(C\) has
  finite colimits which are stable under pullback and preserved by any
  morphism in \(\myFinIndType\).
\end{myProposition}

\begin{myProof}
  Let \((C \myComma T \myComma A) \myElemOf \myFinLimSem\). For a
  finite simplicial type \(I\), let
  \(t \myElemOf I \myPower_{C} \myRepClf_{C} \myMorphism T\) be the
  composite
  \(I \myPower_{C} \myRepClf_{C} \myXMorphism{I \myPower_{C} (a
    \myMapsTo a \myAct_{\myUnivBase_{T}} \myUnivTot_{T})} I
  \myPower_{C} (\myUnivBase_{T} \mySlice_{C} T) \myXMorphism{\myColim}
  (\myUnivBase_{T} \mySlice_{C} T) \myMorphism T\). Observe that
  \(\myApp{t^{\myStar}}{A} \myEquiv I^{\myRCone} \myPower_{C}
  \myRepClf_{C}\), where \(I^{\myRCone}\) is obtained from \(I\) by
  adjoining a terminal vertex. The fiber of
  \(\myApp{t^{\myStar}}{A} \myMorphism I \myPower_{C} \myRepClf_{C}\)
  over \(c \myElemOf C\) is thus the forgetful functor
  \(I^{\myRCone} \myPower (C \mySlice c) \myMorphism I \myPower (C
  \mySlice c)\). Its left adjoint assigns a colimit of \(x\) to each
  diagram \(x \myElemOf I \myMorphism (C \mySlice c)\).
\end{myProof}

For an object \(A \myElemOf C\), the suspension
\(\myApp{\mySuspension}{A}\) is the fiber coproduct of two copies of
\(A \myMorphism \myTerminal\). The \(n\)-sphere is inductively defined
as \(\mySphere^{-1} \myDefEq \myInitial\) and
\(\mySphere^{n \myPlus 1} \myDefEq
\myApp{\mySuspension}{\mySphere^{n}}\).

\begin{myCorollary}
  Let \(C\) be a category with higher inductive types. Then \(C\) has
  suspensions and the \(n\)-sphere for all \(n \myGe -1\) which are
  stable under pullback and preserved by any morphism in
  \(\myFinIndType\). \myQED
\end{myCorollary}

\begin{myDefinition}
  Let \(C\) be a category with finite limits. A \myDefine{natural
    number object} is a tuple
  \((\myNat \myComma \myNatZero \myComma \myNatSuc)\) of
  \(\myNat \myElemOf C\),
  \(\myNatZero \myElemOf \myTerminal \myMorphism \myNat\), and
  \(\myNatSuc \myElemOf \myNat \myMorphism \myNat\) such that for any
  \(x \myElemOf C\), \(z \myElemOf \myTerminal \myMorphism x\), and
  \(s \myElemOf x \myMorphism x\), the type
  \[
    (h \myElemOf \myApp{\myHom}{\myNat \myComma x}) \myBinProd (h
    \myComp \myNatZero \myId z) \myBinProd (h \myComp \myNatSuc \myId
    s \myComp h)
  \]
  is contractible. Equivalently, it is an initial object in the
  category of such tuples \((x \myComma z \myComma s)\).
\end{myDefinition}

\begin{myProposition}
  \label{prop-ind-type-natural-numbers}
  Let \(C\) be a category with higher inductive types. Then \(C\) has
  a natural number object which is stable under pullback and preserved
  by any morphism in \(\myFinIndType\).
\end{myProposition}

\begin{myProof}
  For \((C \myComma T \myComma A) \myElemOf \myFinLimInit\), we define
  \(t \myElemOf C \myMorphism \myApp{\myOp_{C}}{T}\) to be
  \(\myUnivTot_{T} \myBinProd_{\myUnivBase_{T}} (\myUnivTot_{T}
  \myExp_{\myUnivBase_{T}} \myUnivTot_{T})\). Then
  \(\myApp{t^{\myStar}}{A}\) is
  \(\myRepClfPt_{C} \myBinProd_{\myRepClf_{C}} (\myRepClfPt_{C}
  \myExp_{\myRepClf_{C}} \myRepClfPt_{C})\). Its fiber over
  \(c \myElemOf C\) is the category of tuples
  \((x \myComma z \myComma s)\) of \(x \myElemOf (C \mySlice c)\),
  \(z \myElemOf \myTerminal \myMorphism x\), and
  \(s \myElemOf x \myMorphism x\). Thus, the left adjoint of
  \(\myApp{t^{\myStar}}{A} \myMorphism C\) assigns a natural number
  object in \((C \mySlice c)\) to every \(c \myElemOf C\).
\end{myProof}

\begin{myDefinition}
  Let \(C\) be a category with finite limits. We inductively define
  \myDefine{\(n\)-truncatedness} of a morphism in \(C\) for \(n \myGe
  -2\) as follows.
  \begin{itemize}
  \item A morphism is \((-2)\)-truncated if it is an equivalence.
  \item A morphism \(f \myElemOf y \myMorphism x\) is
    \((n \myPlus 1)\)-truncated if the diagonal morphism
    \(y \myMorphism y \myBinProd_{x} y\) is \(n\)-truncated.
  \end{itemize}
  We say an object \(x \myElemOf C\) is \myDefine{\(n\)-truncated} if
  the morphism \(x \myMorphism \myTerminal\) is \(n\)-truncated.
\end{myDefinition}

\begin{myDefinition}
  Let \(C\) be a category with finite limits, let \(x \myElemOf C\),
  and let \(n \myGe -2\). An \myDefine{\(n\)-truncation of \(x\)} is
  an initial object in the full subcategory of \((x \mySlice C)\)
  consisting of those \(x \myMorphism y\) such that \(y\) is
  \(n\)-truncated. The \(n\)-truncation of \(x\) is denoted by
  \(\myTrunc{n}{x}\) if it exists.
\end{myDefinition}

\begin{myProposition}
  \label{prop-ind-type-truncation}
  Let \(C\) be a category with higher inductive types. Then \(C\) has
  \(n\)-truncations which are stable under pullback and preserved by
  any morphism in \(\myFinIndType\).
\end{myProposition}

\begin{myProof}
  Let \((C \myComma T \myComma A) \myElemOf \myFinLimInit\). For
  \(n \myGe -2\), we first construct an object
  \(\myTruncatedType_{n} \myElemOf (\myApp{\myOp_{C}}{T} \mySlice_{C}
  \myUnivBase_{T})\) as follows. \(\myTruncatedType_{-2}\) is
  \(\myIntern{\myTerminal} \myElemOf \myTerminal \myMorphism
  \myUnivBase_{T}\). We define \(\myTruncatedType_{n \myPlus 1}\) to
  be
  \(\myApp{p_{\myStar}}{\myApp{\myIntern{\myId}^{\myStar}}{\myTruncatedType_{n}}}\)
  where
  \(p \myElemOf \myUnivTot_{T} \myBinProd_{\myUnivBase_{T}}
  \myUnivTot_{T} \myMorphism \myUnivBase_{T}\) is the projection which
  is exponentiable. Every \(\myTruncatedType_{n}\) is a subobject of
  \(\myUnivBase_{T}\), and for \(x \myElemOf \myApp{\myOp_{C}}{T}\), a
  morphism \(f \myElemOf x \myMorphism \myUnivBase_{T}\) factors
  through \(\myTruncatedType_{n}\) if and only if the morphism
  \(\myApp{f^{\myStar}}{\myUnivTot_{T}} \myMorphism x\) is
  \(n\)-truncated. We then define
  \(t_{n} \myElemOf \myRepClf_{C} \myMorphism T\) by the fiber
  coproduct
  \[
    \begin{tikzcd}
      \myUnivBase_{T}
      \arrow[r]
      \arrow[d]
      & \myTruncatedType_{n}
      \arrow[d]
      \\ a \myAct_{\myUnivBase_{T}} \myUnivTot_{T}
      \arrow[r]
      & \myApp{t_{n}}{a}.
    \end{tikzcd}
  \]
  The fiber of \(\myApp{t_{n}^{\myStar}}{A}\) over \(c \myElemOf C\)
  is the category of tuples \((a \myComma x \myComma f)\) of
  \(a \myComma x \myElemOf (C \mySlice c)\) and
  \(f \myElemOf a \myMorphism x\) such that \(x\) is
  \(n\)-truncated. The left adjoint of
  \(\myApp{t_{n}^{\myStar}}{A} \myMorphism \myRepClf_{C}\) thus
  assigns an \(n\)-truncation of \(a\) to each
  \(a \myElemOf (C \mySlice c)\).
\end{myProof}

\((-1)\)-truncated objects are called propositions, and
\((-1)\)-truncation is called propositional truncation. For
propositions \(P\) and \(Q\) in \(C\), the disjunction
\(P \myLogicOr Q\) is defined as \(\myTrunc{-1}{P \myBinCoprod
  Q}\). For an object \(A \myElemOf C\) and a proposition
\(P \myElemOf (C \mySlice A)\), the existential quantification
\(\myApp{\myExists_{A}}{P}\) is defined as
\(\myTrunc{-1}{\myApp{A_{\myBang}}{P}}\), where
\(A_{\myBang} \myElemOf (C \mySlice A) \myMorphism C\) is the
forgetful functor.

\begin{myCorollary}
  Let \(C\) be a category with higher inductive types. Then \(C\) has
  disjunctions and existential quantifications of propositions which
  are stable under pullback and preserved by any morphism in
  \(\myFinIndType\). \myQED
\end{myCorollary}

\(0\)-truncated objects are called sets, and \(0\)-truncation is
called set truncation. For \(n \myGe 0\), the \(n\)-th loop space
\(\myApp{\myLoopSpace^{n}}{A \myComma a}\) of a pointed object
\((A \myElemOf C \myComma a \myElemOf \myTerminal \myMorphism A)\) is
inductively defined as
\(\myApp{\myLoopSpace^{0}}{A \myComma a} \myDefEq A\) and
\(\myApp{\myLoopSpace^{n \myPlus 1}}{A \myComma a} \myDefEq
\myApp{\myLoopSpace^{n}}{a \myId a \myComma \myRefl}\), where
\(a \myId a\) is the fiber product of two copies of \(a\) and
\(\myRefl \myElemOf \myTerminal \myMorphism a \myId a\) is the
diagonal morphism. The \(n\)-th homotopy group
\(\myApp{\myHoGrp_{n}}{A \myComma a}\) of \((A \myComma a)\) is
defined as \(\myTrunc{0}{\myApp{\myLoopSpace^{n}}{A \myComma
    a}}\). Suspensions \(\myApp{\mySuspension}{A}\) and in particular
the \(k\)-sphere \(\mySphere^{k}\) for \(k \myGe 0\) are regarded as
pointed objects with the first inclusion
\(\myTerminal \myMorphism \myApp{\mySuspension}{A}\).

\begin{myCorollary}
  Let \(C\) be a category with higher inductive types. Then \(C\) has
  \(n\)-th homotopy groups of pointed objects, in particular of the
  \(k\)-sphere for \(k \myGe 0\), which are stable under pullback and
  preserved by any morphism in \(\myFinIndType\). \myQED
\end{myCorollary}

\section{Canonicity for higher inductive types}
\label{sec-canonicity}

We prove canonicity for higher inductive types
(\myRef{prop-inductive-type-canonicity}). We follow the standard
method using the so-called Freyd cover or Sierpinski cone
\myCite{lambek1986introduction,shulman2015inverse}, which is the comma
category \((\myType \mySlice \myGSec)\) for the global section functor
\(\myGSec \myElemOf C \myMorphism \myType\).

\begin{myProposition}
  \(\myType\) has higher inductive types.
\end{myProposition}

\begin{myProof}
  By \myRef{prop-when-lex-has-ind-type}, it suffices to construct a
  limit metatheory \(T\) over \(\myType\) and a semantics \(A\) of
  \(T\) that admits initial algebras. Let \(A' \myMorphism T'\) be the
  unstraightening of the semantics functor
  \(\myFinLimThAlg \myElemOf \myApp{\myOp}{\myFinLimTh} \myMorphism
  \myEnlarge \myCat\) introduced in \myRef{sec-limit-theories}. Since
  \(\myType\) is the completion of the point under small colimits,
  \(\myEnlarge \myCat\) is regarded as the full subcategory of
  \(\myApp{\myOp}{\myType} \myPower \myEnlarge \myCat\) consisting of
  those functors preserving small limits. The image of
  \(A' \myMorphism T'\) by the inclusion
  \(\myEnlarge \myCat \myMorphism \myApp{\myOp}{\myType} \myPower
  \myEnlarge \myCat \myEquiv \myFiber{(\myEnlarge \myCart)}{\myType}\)
  is what we need.
\end{myProof}

\begin{myLemma}
  \label{prop-ind-type-mor-retract}
  Let \(C\) and \(D\) be categories with higher inductive types. Then
  the class of functors \(C \myMorphism D\) in \(\myFinIndType\) is
  closed in \(C \myPower D\) under retract.
\end{myLemma}

\begin{myProof}
  Let \(F \myElemOf C \myMorphism D\) be a morphism in
  \(\myFinIndType\) and let \(F'\) be a retract of \(F\) in
  \(C \myPower D\). We show that \(F'\) is in \(\myLex\) and
  \(\myApp{\myFinLimThAlg}{F'}\) is in \(\myFinLimInit\). By
  \myRef{prop-adj-locally-closed-under-retract}, we see that
  \(\myLex \mySub \myCat\) and \(\myFinLimInit \mySub \myFinLimSem\)
  are locally full and locally closed under retract, so \(F'\) is in
  \(\myLex\). As in \myRef{cst-2pr-left-2-functor}, the morphism part
  of \(\myFinLimThAlg \myElemOf \myLex \myMorphism \myFinLimSem\) is
  extended to a functor. Thus, \(\myApp{\myFinLimThAlg}{F'}\) is a
  retract of \(\myApp{\myFinLimThAlg}{F}\) in
  \(\myApp{\myHHom_{\myFinLimSem}}{C \myComma D}\) and thus belongs to
  \(\myFinLimInit\).
\end{myProof}

\begin{myLemma}
  \label{prop-comma-left-adjoint}
  Let
  \[
    \begin{tikzcd}
      C
      \arrow[r, "F"]
      \arrow[d, "P"']
      & C'
      \arrow[d, "P'"]
      \\ D
      \arrow[r, "G"']
      & D'
    \end{tikzcd}
  \]
  be a commutative square of categories. Let
  \(R \myElemOf (C' \mySlice F) \myMorphism (D' \mySlice G)\) denote
  the induced functor, that is,
  \(\myApp{R}{f \myElemOf x' \myMorphism \myApp{F}{x}} \myDefEq
  (\myApp{P'}{x'} \myXMorphism{\myApp{P'}{f}} \myApp{P'}{\myApp{F}{x}}
  \myEquiv \myApp{G}{\myApp{P}{x}})\). Suppose that \(P\) and \(P'\)
  have left adjoints. Then \(R\) has a left adjoint, and the mates of
  the projection squares from \(R\) to \(P'\) and from \(R\) to \(P\)
  are equivalences.
\end{myLemma}

\begin{myProof}
  Let \(Q\) and \(Q'\) be the left adjoints of \(P\) and \(P'\),
  respectively. Let
  \(\alpha \myElemOf Q' \myComp G \myMMorphism F \myComp Q\) denote the
  mate. The left adjoint
  \(S \myElemOf (D' \mySlice G) \myMorphism (C' \mySlice F)\) of \(R\)
  is then defined by
  \(\myApp{S}{g \myElemOf y' \myMorphism \myApp{G}{y}} \myDefEq
  (\myApp{Q'}{y'} \myXMorphism{\myApp{Q'}{g}} \myApp{Q'}{\myApp{G}{y}}
  \myXMorphism{\alpha} \myApp{F}{\myApp{Q}{y}})\). By construction, the
  mates of the two projection squares are equivalences.
\end{myProof}

\begin{myLemma}
  \label{prop-ind-type-comma}
  Let \(C\) and \(D\) be categories with higher inductive types and
  let \(F \myElemOf C \myMorphism D\) be a functor preserving finite
  limits. Then \((D \mySlice F)\) has higher inductive types, and the
  projections \((D \mySlice F) \myMorphism D\) and
  \((D \mySlice F) \myMorphism C\) are in \(\myFinIndType\).
\end{myLemma}

\begin{myProof}
  Let
  \(((D \mySlice F) \myComma T \myComma A) \myElemOf \myFinLimSem\) be
  \((\myApp{\myFinLimThAlg}{D} \mySlice
  \myApp{\myFinLimThAlg}{F})\). Because this comma object in
  \(\myFinLimSem\) is computed component-wise,
  \myRef{prop-comma-left-adjoint} implies that \(A \myMorphism T\) has
  a left adjoint \(H\) computed component-wise. It then follows that
  \(H\) is cartesian over \((D \mySlice C)\) and the unit is an
  equivalence. Therefore, \(A\) has initial algebras. Then
  \((D \mySlice F)\) has higher inductive types by
  \myRef{prop-when-lex-has-ind-type}. By construction, the two
  projections are in \(\myFinIndType\).
\end{myProof}

\begin{myTheorem}
  \label{prop-inductive-type-canonicity}
  Let \(C\) be the initial object in \(\myFinIndType\). Then the
  global section functor \(\myGSec \myElemOf C \myMorphism \myType\)
  defined by
  \(\myApp{\myGSec}{x} \myDefEq \myApp{\myHom}{\myTerminal \myComma
    x}\) is in \(\myEnlarge \myFinIndType\).
\end{myTheorem}

\begin{myProof}
  By \myRef{prop-ind-type-comma}, \((\myType \mySlice F)\) has higher
  inductive types, and the projections
  \(\myProj_{1} \myElemOf (\myType \mySlice F) \myMorphism \myType\)
  and \(\myProj_{2} \myElemOf (\myType \mySlice F) \myMorphism C\) are
  in \(\myEnlarge \myFinIndType\). By the initiality of \(C\), the
  projection \(\myProj_{2}\) has a unique section \(G'\) in
  \(\myEnlarge \myFinIndType\). Let
  \(G \myDefEq \myProj_{1} \myComp G'\). We have a natural
  transformation \(\epsilon \myElemOf G \myMMorphism \myGSec\). Recall that
  \(\myGSec\) is the initial object in the full subcategory of
  \(C \myPower \myType\) consisting of those functors preserving
  finite limits. Then \(\epsilon\) has a unique section, and thus
  \(\myGSec\) is a retract of \(G\). By
  \myRef{prop-ind-type-mor-retract}, \(\myGSec\) is in
  \(\myEnlarge \myFinIndType\).
\end{myProof}

For a proposition \(P\) in \(C\), we say \(P\) is inhabited if there
exists a morphism \(\myTerminal \myMorphism P\).

\begin{myCorollary}
  \label{prop-classical-canonicity}
  Let \(C\) be the initial object in \(\myFinIndType\). The following
  hold.
  \begin{description}
  \item[Canonicity for the natural number object] For every morphism
    \(f \myElemOf \myTerminal \myMorphism \myNat\) in \(C\), there
    exists a unique natural number \(n\) such that
    \(f \myId \myNatSuc^{n} \myComp \myNatZero\).
  \item[Disjunction property] Let \(P \myComma Q \myElemOf C\) be
    propositions. If \(P \myLogicOr Q\) is inhabited, then either
    \(P\) or \(Q\) is inhabited.
  \item[Existence property] Let \(A \myElemOf C\) and let
    \(P \myElemOf (C \mySlice A)\) be a proposition. If
    \(\myApp{\myExists_{A}}{P}\) is inhabited, then there exists a
    morphism \(a \myElemOf \myTerminal \myMorphism A\) such that
    \(\myApp{a^{\myStar}}{P}\) is inhabited.  \myQED
  \end{description}
\end{myCorollary}

\begin{myCorollary}
  \label{prop-homotopy-group-canonicity}
  Let \(C\) be the initial object in \(\myFinIndType\). For any
  \(n \myGe 0\) and \(k \myGe 0\), we have a canonical equivalence
  \(\myApp{\myHoGrp_{n}}{\mySphere^{k}} \myEquiv
  \myApp{\myHom_{C}}{\myTerminal,
    \myApp{\myHoGrp_{n}}{\mySphere^{k}}}\) (of sets for \(n \myId 0\)
  and of groups for \(n \myGe 1\)). \myQED
\end{myCorollary}

The proof of \myRef{prop-inductive-type-canonicity} can immediately be
generalized.

\begin{myTheorem}
  \label{prop-inductive-type-canonicity-general}
  Let \(X \mySub \myLex\) be a subobject in \(\myPrRFin\). Let \(C\)
  be the initial object in \(X \myBinProd_{\myLex} \myFinIndType\) and
  let \(\myGSec \myElemOf C \myMorphism \myType\) denote the global
  section functor. Suppose that the comma category
  \((\myType \mySlice \myGSec)\) and the projection
  \((\myType \mySlice \myGSec) \myMorphism C\) are in
  \(\myEnlarge X\). Then \(\myGSec\) is in
  \(\myEnlarge \myFinIndType\). \myQED
\end{myTheorem}

\begin{myDefinition}
  A category \(C\) with finite limits is \myDefine{locally cartesian
    closed} if every
  \(f \myElemOf \myApp{\myObj}{\myCell_{1} \myPower C}\) is an
  exponentiable morphism in \(C\). We have the
  \(\myFinitary\)-presentable category \(\myLCCC\) of locally
  cartesian closed categories.
\end{myDefinition}

\begin{myRemark}
  We would not say ``every \(f \myElemOf \myCell_{1} \myPower C\) is
  an exponentiable morphism in \(C\)''. This reflects the fact that
  \((f \myElemOf y \myMorphism x \myComma g \myElemOf z \myMorphism y)
  \myMapsTo \myApp{f_{\myStar}}{g}\) is not functorial. Because of
  this, \(\myLCCC\) is defined as just a category, and there is no way
  to make it a well-behaved \(2\)-category.
\end{myRemark}

\begin{myTheorem}
  \label{prop-inductive-type-canonicity-lccc}
  Let \(C\) be the initial object in
  \(\myLCCC \myBinProd_{\myLex} \myFinIndType\). Then the global
  section functor \(\myGSec \myElemOf C \myMorphism \myType\) is in
  \(\myEnlarge \myFinIndType\).
\end{myTheorem}

\begin{myProof}
  By \myRef{prop-inductive-type-canonicity-general}, this follows from
  the fact that every morphism in \((\myType \mySlice \myGSec)\) is
  exponentiable and the projection
  \((\myType \mySlice \myGSec) \myMorphism C\) preserves
  exponentials. See, for example,
  \myCite{shulman2015inverse,uemura2017fibred} for the computation of
  exponentials in \((\myType \mySlice \myGSec)\).
\end{myProof}

\printbibliography

\end{document}